\documentclass{article}
\usepackage[T1]{fontenc}
\usepackage{geometry}
\usepackage{amssymb}
\usepackage{amsmath,cases}
\usepackage{amsthm}
\usepackage{enumitem}
\usepackage{comment}
\usepackage{xcolor}
\usepackage{graphicx}
\usepackage{textcomp}
\usepackage{mathrsfs}
\usepackage{bbm}
\usepackage{bm}
\usepackage{caption}
\usepackage{stmaryrd}
\usepackage{url}
\usepackage{calc}
\usepackage{ifthen}
\usepackage{psfrag}
\usepackage[active]{srcltx}   
\usepackage{pdfsync}
\usepackage[colorlinks=true,bookmarksnumbered=true, pdfauthor={J\'er\'emie Bettinelli}]{hyperref}									
\usepackage{todonotes}

\newcommand{\lv}[1]{{#1}}
\newcommand{\sv}[1]{}

\newcommand{\de}{\mathrel{\mathop:}\hspace*{-.6pt}=}
\newcommand{\ooo}{\bullet}
\newcommand{\dt}{\rightarrowtriangle}
\newcommand{\da}{\leftarrowtriangle}

\newcommand{\ba}{\bm{a}}

\newcommand{\cM}{\mathcal{M}}

\newcommand{\m}{\mathfrak{m}}

\newcommand{\pc}{\mathbbm c}
\newcommand{\ph}{\mathbbm h}
\newcommand{\pl}{\mathbbm l}
\newcommand{\pp}{\mathbbm p}
\newcommand{\pr}{\mathbbm r}
\newcommand{\rev}{\operatorname{rev}}

\definecolor{vert}{rgb}{0,0.666,0}
\definecolor{violet}{rgb}{0.8,0,0.8}

\theoremstyle{plain}
\newtheorem{thm}{Theorem}
\newtheorem{prop}[thm]{Proposition}
\newtheorem{defi}{Definition}
\newtheorem*{wrn}{Warning}

\theoremstyle{definition}
\newtheorem{rem}{Remark}
\newtheorem*{note}{Note}
\newtheorem*{ack}{Acknowledgment}

\newenvironment{pre}[1][\proofname]{%
  \proof[#1]%
}{\endproof}

\makeatletter
\renewcommand*{\@fnsymbol}[1]{\ensuremath{\ifcase#1\or \dagger\or \ddagger\or
   \mathsection\or \mathparagraph\or \|\or **\or \dagger\dagger
   \or \ddagger\ddagger \else\@ctrerr\fi}}
\makeatother
\title{Slit-slide-sew bijections for bipartite and quasibipartite plane maps}
\author{J\'er\'emie Bettinelli\thanks{cnrs \& Laboratoire d'Informatique de l'\'Ecole polytechnique; \href{mailto:jeremie.bettinelli@normalesup.org}{\nolinkurl{jeremie.bettinelli@normalesup.org}};\newline \nolinkurl{www.normalesup.org/}\texttildelow\nolinkurl{bettinel}. This work is partially supported by Grant ANR-14-CE25-0014 (GRAAL).}}

\includecomment{lvc}
\excludecomment{svc}

\begin{document}
\maketitle

\begin{abstract}
We unify and extend previous bijections on plane quadrangulations to bipartite and quasibipartite plane maps. Starting from a bipartite plane map with a distinguished edge and two distinguished corners (in the same face or in two different faces), we build a new plane map with a distinguished vertex and two distinguished half-edges directed toward the vertex. The faces of the new map have the same degree as those of the original map, except at the locations of the distinguished corners, where each receives an extra degree. The idea behind this bijection is to build a path from the distinguished elements, slit the map along it, and sew back after sliding by one unit, thus mildly modifying the structure of the map at the extremities of the sliding path. This bijection provides a sampling algorithm for uniform maps with prescribed face degrees and allow to recover Tutte's famous counting formula for bipartite and quasibipartite plane maps.

In addition, we explain how to decompose the previous bijection into two more elementary ones, which each transfer a degree from one face of the map to another face. In particular, these transfer bijections are simpler to manipulate than the previous one and this point of view simplifies the proofs.
\end{abstract}

\newcommand{\noteM}[1]{{\bfseries\textcolor{blue}{~(M: #1)}}}

\section{Introduction}

This paper is a sequel to~\cite{bettinelli14ifq}, in which we presented two bijections on plane quadrangulations with a boundary. In the present work, we show how to generalize these bijections to bipartite and, in some cases, quasibipartite plane maps. Recall that a plane map is an embedding of a finite connected graph (possibly with multiple edges and loops) into the sphere, considered up to orientation-preserving homeomorphisms. It is \emph{bipartite} if every of its faces have an even degree and \emph{quasibipartite} if it has two faces of odd degree and all other faces of even degree. \lv{Note that, as the sum of the face degrees equals twice the number of edges, the number of faces with an odd degree must be even, so that quasibipartite maps are the simplest maps to consider after bipartite maps.}

\begin{figure}[ht]
		\psfrag{f}[][][.9]{$f_1$}
		\psfrag{a}[][][.9]{$f_{10}$}
		\psfrag{z}[][][.9]{$f_8$}
		\psfrag{e}[][][.9]{$f_3$}
		\psfrag{r}[][][.9]{$f_{12}$}
		\psfrag{t}[][][.9]{$f_{14}$}
		\psfrag{y}[][][.9]{$f_{11}$}
		\psfrag{u}[][][.9]{$f_9$}
		\psfrag{i}[][][.9]{$f_2$}
		\psfrag{o}[][][.9]{$f_{13}$}
		\psfrag{p}[][][.9]{$f_7$}
		\psfrag{q}[][][.9]{$f_5$}
		\psfrag{s}[][][.9]{$f_4$}
		\psfrag{d}[][][.9]{$f_6$}
	\centering\includegraphics[width=9cm]{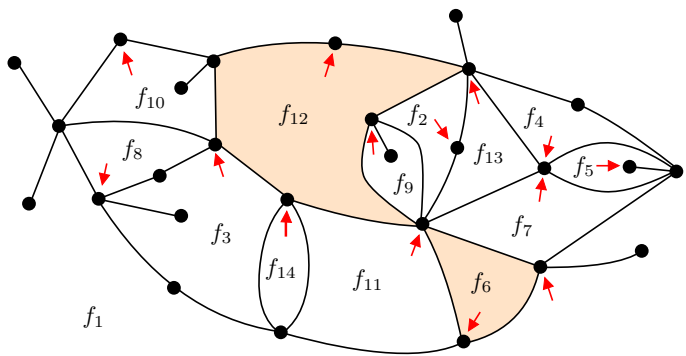}
	\caption{A \sv{quasibipartite }map of type $(20,4,8,4,4,3,4,4,4,6,4,7,4,2)$.\lv{ It is quasibipartite because it has exactly two faces of odd degree ($f_6$ and~$f_{12}$); throughout the paper, we will highlight odd-degree faces by coloring them orange. Each face has a marked corner (represented by a red arrowhead).}}
	\label{maptypea}
\end{figure}

The number of such maps with prescribed face degrees has been computed by several methods. For an $r$-tuple $ \ba=( a_1,\ldots, a_r)$ of positive integers, let us denote by $M(\ba)$ the number of plane maps with~$r$ numbered faces~$f_1$, \ldots, $f_r$ of respective degrees~$ a_1$, \ldots, $ a_r$, where each face has a marked corner. The $r$-tuple~$\ba$ will be called the \emph{type} of such maps (see Figure~\ref{maptypea}). By elementary considerations and Euler's characteristic formula, the integers
\begin{equation}\label{ev}
E(\ba)\de\frac12 \sum_{i=1}^r a_i\qquad\text{ and }\qquad V(\ba)\de E(\ba)-r+2
\end{equation}
are respectively the numbers of edges and vertices of maps of type~$\ba$. Solving a technically involved recurrence, Tutte~\cite{tutte62cs} showed that, when at most two~$a_i$'s are odd, that is, for bipartite or quasibipartite maps, 
\begin{equation}\label{slicings}
M(\ba)=\frac{\big(E(\ba)-1\big)!}{V(\ba)\,!}\prod_{i=1}^r\alpha( a_i),\qquad\text{ where }\quad\alpha(x)\de\frac{x!}{\big\lfloor x/2\big\rfloor!\,\big\lfloor(x-1)/2\big\rfloor!}\,.
\end{equation}
Formula~\eqref{slicings}, commonly referred to as \emph{Tutte's formula of slicings}, was later recovered by Cori~\cite{cori75code,cori76} thanks to a so-called \emph{transfer bijection}, roughly consisting in iteratively transferring one degree from a face to a neighboring face, until the map has a very simple structure. Using a bijective encoding by so-called \emph{blossoming trees}, Schaeffer~\cite{schaeffer97bij} then recovered it in the bipartite case. Finally, we may also obtain it by using the so-called \emph{Bouttier--Di Francesco--Guitter bijection} \cite{bouttier04pml}, which encodes plane maps by tree-like structures called \emph{mobiles}: see~\cite{CoFu14} for the computation of related generating functions using this approach.

In the present work, we give a bijective interpretation for the following combinatorial identities, which somehow allows to ``grow'' maps by adding to a bipartite map two new corners either to the same face or to two different faces. 

\begin{prop}[Adding two corners to the same face]
\label{propp2}
Let $\ba=(a_1,\ldots, a_r)$ be an $r$-tuple of positive even integers and let $\tilde \ba=(\tilde a_1,\ldots,\tilde a_r)\de( a_1+2, a_2,\ldots, a_r)$. Then the following identity holds:
\begin{align}\label{eqp2}
(a_1+1)\,(a_1+2)\, E(\ba)\, M(\ba)&=\big\lfloor {\tilde a_1}/{2}\big\rfloor\big\lfloor{(\tilde a_1-1)}/{2}\big\rfloor\, V(\tilde \ba)\,M(\tilde \ba).
\end{align}
\end{prop}

\begin{prop}[Adding one corner to each of two different faces]
\label{propp1p1}
Let $\ba=(a_1,\ldots, a_r)$ be an $r$-tuple of positive even integers and let $\tilde \ba=(\tilde a_1,\ldots,\tilde a_r)\de( a_1+1, a_2+1,a_3,\ldots, a_r)$. Then the following identity holds:
\begin{align}\label{eqp1p1}
(a_1+1)\,(a_2+1)\, E(\ba)\, M(\ba)&=\big\lfloor {\tilde a_1}/{2}\big\rfloor\big\lfloor{\tilde a_2}/{2}\big\rfloor\, V(\tilde \ba)\,M(\tilde \ba).
\end{align}
\end{prop}

For the $r$-tuple $(2,\ldots,2)$, it is easy to see that $M(2,\ldots,2)=2^{r-1} (r-1)!$ as there is only one map with~$r$ faces of degree~$2$ and a chosen first face with its marked corner, and there are $(r-1)!$ ways to order the remaining faces and $2^{r-1}$ ways to choose the remaining marked corners. This initial condition, together with the above propositions and the obvious exchangeability of the coordinates of~$\ba$ provides yet another proof of~\eqref{slicings}.

\bigskip
We will use the technique introduced in~\cite{bettinelli14ifq} of what we call \emph{slit-slide-sew bijections}, and whose idea is the following. We will interpret the sides of~\eqref{eqp2} and~\eqref{eqp1p1} as counting maps with some distinguished ``elements.'' More precisely, in each case, the term in~$M$ counts maps of some type and the three prefactors will count something whose number only depends on this type: it can be a corner, an edge, a vertex, or something a bit more intricate. For instance, the left-hand side of~\eqref{eqp1p1} counts maps of type~$\ba$ with a distinguished corner in~$f_1$, a distinguished corner in~$f_2$ and a distinguished edge (for any~$i$, there are $a_i+1$ corners in~$f_i$ because of the already marked corner; see Section~\ref{secprel} for the convention on distinguishing corners).

From a map with its distinguished elements, we first construct a directed path. We then slit the map along this path and we sew back together the sides of the slit but after sliding by one unit. Let us look at a face lying to the left of some edge of the path. Before the operation, it is adjacent to the face lying to the right of the same edge and, after the operation, it is adjacent to the face lying to the right of the next or previous edge along the path. This operation mildly modifies the map along the path but does not affect its faces, except around the extremities of the path. In the process, new distinguished
elements naturally appear in the resulting map. Plainly, in order for this operation to work, the path we construct has to be totally recoverable
from the new distinguished elements.

\bigskip

We will furthermore see the previous bijections as compositions of two more elementary bijections, which can be thought of as ``transferring'' a corner from a face, say~$f_{r+1}$, to another face, say~$f_1$. In the case where~$f_{r+1}$ has degree~$1$, it somehow vanishes into a vertex. We chose to use an $r+1$-th face for these operations as we will see the previous mappings as compositions of the following ones by using an extra face. More precisely, by a slight modification, we may transform a distinguished edge into an extra degree-2 face and use twice the bijections interpreting the following identities in order to transfer both corners of the extra face to the desired faces.

\begin{prop}[Transferring a corner from a face of degree at least~$2$]
\label{propp1m1}
Let $\ba=(a_1,\ldots, a_{r+1})$ be an $r+1$-tuple of positive integers with~$a_{r+1}\ge 2$, and either all even or such that only~$a_{r+1}$ and one other coordinate are odd. Let also $\tilde \ba=(\tilde a_1,\ldots,\tilde a_{r+1})\de( a_1+1, a_2,\ldots, a_{r}, a_{r+1}-1)$. Then the following identity holds:
\begin{align}\label{eqp1m1}
(a_1+1)\,\big\lfloor {a_{r+1}}/{2}\big\rfloor\, M(\ba)&=\big\lfloor{\tilde a_1}/{2}\big\rfloor\,(\tilde a_{r+1}+1)\, M(\tilde \ba).
\end{align}
\end{prop}

\begin{prop}[Transferring a corner from a degree~$1$-face]
\label{propp1m10}
Let $\ba=(a_1,\ldots, a_r,1)$ be an $r+1$-tuple of positive integers with two odd coordinates and let $\tilde \ba=(\tilde a_1,\ldots,\tilde a_r)\de( a_1+1, a_2,\ldots, a_r)$. Then the following identity holds:
\begin{align}\label{eqp1m10}
(a_1+1)\, M(\ba)&=\big\lfloor{\tilde a_1}/{2}\big\rfloor\,V(\tilde \ba)\,M(\tilde \ba).
\end{align}
\end{prop}

\lv{
The left-hand side of~\eqref{eqp1m10} may seem to miss a factor but really, one should see the $r+1$-th face as a second distinguished element, so that there always are two distinguished elements in both sides of~\eqref{eqp1m1} and~\eqref{eqp1m10}.
}

\paragraph{Related works.}
Let us mention at this point that our bijections bear some similarities with two related works. In the papers we mentioned earlier, Cori~\cite{cori75code,cori76} also transfers one degree from a face to another one. In his approach, he does so in a local way, in the sense that the degree passes from a face to one of its neighbor. In the present work, our transfer bijections are global in the sense that the degree passes from a face to an arbitrarily far away one. Moreover, the notion of geodesic path along which we slide the map is of crucial importance.

In a very recent work, Louf~\cite{louf18bij} introduced a new family of bijections accounting for formulas on plane maps arising from the so-called KP hierarchy. His bijections also strongly rely on the mechanism of sliding along a path but, in his case, the path is also somehow local (although arbitrary long) as it is canonically defined from only one vertex using a depth-first search exploration of the map. Another difference of importance is that his mappings may produce two maps as an output, which corresponds to the fact that the formulas in question are quadratic; in the present work, the output is always one map, which corresponds to linear formulas.

\lv{
\paragraph{Structure.}
The remaining of the paper is organized as follows. We start by recalling in Section~\ref{secprel} the definitions and conventions we use, as well as some elementary facts on bipartite and quasibipartite plane maps. We will next see in Section~\ref{secbb} and~\ref{secbq} the bijections that account for Propositions~\ref{propp2} and~\ref{propp1p1}. In Section~\ref{sectrans}, we present the transfer bijections interpreting Propositions~\ref{propp1m1} and~\ref{propp1m10} and explain how our previous bijections can be decomposed as two such bijections. We explain in Section~\ref{secsample} how to sample a uniform map of a given type using our bijections. Finally, Section~\ref{secopen} is devoted to the generalization of our bijections.

\begin{wrn}
Throughout the paper, we will present several bijections between sets of maps carrying distinguished elements. In order to lighten the notation, we will always denote by~$\cM$ and~$\tilde\cM$ the sets in bijections. The definitions of these sets depend on the section they appear in; they are always clearly defined (with helping pictographs) at the beginning of the section in question.
\end{wrn}

\lv{\begin{ack}
We thank Olivier Bernardi and \'Eric Fusy for interesting discussions on this topic.
\end{ack}}
}

\section{Preliminaries}\label{secprel}

\begin{figure}[ht]
		\psfrag{c}[][]{\textcolor{blue}{$c$}}
		\psfrag{d}[][]{\textcolor{vert}{$c'$}}
		\psfrag{o}[][]{or}
	\centering\includegraphics{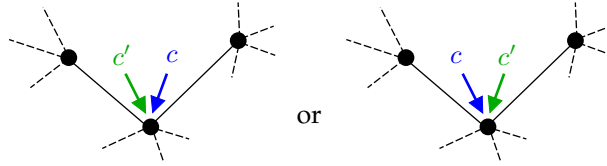}
	\caption{The two different ways of distinguishing twice the same corner.}
	\label{discor}
\end{figure}

We will use the following terminology. We call \emph{half-edge} an edge given with one of its two possible orientations. For a half-edge~$h$, we denote by~$h^-$ its origin, by~$h^+$ its end, and by~$\rev(h)$ its reverse. We say that a half-edge~$h$ is \emph{incident} to a face~$f$ if~$h$ lies on the boundary of~$f$ and has~$f$ to its left. It will be convenient to view corners as half-edges having no origin, only an end. In particular, if~$c$ is a corner, we will write~$c^+$ the vertex corresponding to it\lv{, that is, if~$c$ is the corner delimited by the consecutive half-edges~$h$ and~$h'$, then $c^+\de h^+=h'^-$}. Moreover, we use the convention that distinguishing a corner ``splits'' it into two new corners\lv{. In other words, when we distinguish the same corner for the second time, we have to specify which of its two sides is distinguished}: see Figure~\ref{discor}.

\begin{defi}
A \emph{path} from a vertex~$v$ to a vertex~$v'$ is a finite sequence $\pp=(\pp_1,\pp_2,\dots,\pp_\ell)$ of half-edges such that $\pp_1^-=v$, for $1\le k \le \ell-1$, $\pp_k^+=\pp_{k+1}^-$, and $\pp_\ell^+=v'$. Its \emph{length} is the integer $[\pp]\de\ell$. By convention, the empty path has length~$0$.

A path~$\pp$ is called \emph{self-avoiding} if it does not meet twice the same vertex\sv{.}\lv{, that is, $$\big|\{\pp_1^-,\dots,\pp_{[\pp]}^-,\pp_{[\pp]}^+\}\big|=[\pp]+1.$$

A path~$\pp$ is called a \emph{cycle} if $\pp_{[\pp]}^+=\pp_1^-$.}

The \emph{reverse} of\lv{ a path}~$\pp=(\pp_1,\pp_2,\dots,\pp_\ell)$ is\lv{ the path} $\rev(\pp) \de (\rev(\pp_\ell),\rev(\pp_{\ell-1}),\dots,\rev(\pp_1))$. 
\end{defi}

Let~$\pp$ be a path. We denote by $\pp_{i\to j}$ the path $(\pp_i,\dots,\pp_j)$ if $1\le i\le j\le [\pp]$, or the empty path otherwise. If~$\mathbbm q$ is another path satisfying $\mathbbm q_1^-=\pp_{[\pp]}^+$, we set
$$\pp \ooo \mathbbm q \de (\pp_1,\dots,\pp_{[\pp]},\mathbbm q_1,\dots,\mathbbm q_{[\mathbbm q]})$$
the concatenation of~$\pp$ and~$\mathbbm q$. Throughout this paper, the notion of metric we use is the graph metric: if~$\m$ is a map, the distance $d_\m(v,v')$ between two vertices~$v$ and~$v'$ is the smaller~$\ell$ for which there exists a path of length~$\ell$ from~$v$ to~$v'$. A \emph{geodesic} from~$v$ to~$v'$ is such a path. The \emph{leftmost geodesic} from a half-edge~$h$ (or a corner) to a vertex or to a corner is constructed as follows. First, we consider all the geodesics from~$h^+$ to the vertex or to the vertex corresponding to the corner. We take the set of all the first steps of these geodesics. Starting from~$h$, we select the first half-edge to its left that belongs to this set.\lv{ In other words, we turn clockwise around $h^+$ and select the first half-edge of this set that we meet. Note that this half-edge may be $\rev(h)$ if this is the only half-edge in the set.} Then we iterate the process from this half-edge until we reach the desired vertex.\lv{ Remark that this path may be empty if~$h^+$ is the desired vertex and that it is a geodesic.} The \emph{rightmost geodesic}\lv{ from a half-edge or a corner to a vertex or a corner} is defined in a similar way\lv{, by replacing the word ``left'' with the word ``right'' in the previous definition}.

For two corners~$c$ and~$c'$ and a self-avoiding path~$\pp$ from~$c^+$ to~$c'^+$ in a map~$\m$, we may slit the map~$\m$ along~$\pp$ from~$c$ to~$c'$ by doubling each edge of~$\pp$. In the resulting object, there are two copies of the initial path~$\pp$, one lying to the left of~$\pp$ and one lying to its right. These are respectively called the \emph{left copy} and \emph{right copy} of~$\pp$. See Figure~\ref{bb}. 
\lv{We will intensively use this operation, and even generalize it in due course for slightly more complicated paths.

\begin{rem}\label{remdiscon}
Note that, if~$c$ and~$c'$ do not lie in the same face of~$\m$, the resulting object is a map, whereas it consists in two separate maps if~$c$ and~$c'$ lie in the same face of~$\m$.
\end{rem}
}

We say that a half-edge~$h$ is \emph{directed toward} a vertex~$v$ if $d_\m(h^+,v)<d_\m(h^-,v)$, that it is \emph{directed away from}~$v$ if $d_\m(h^+,v)>d_\m(h^-,v)$ and that it is \emph{parallel} to~$v$ if $d_\m(h^+,v)=d_\m(h^-,v)$. In the following figures and pictographs, we will represent half-edges with half arrowheads and use the shorthand notation $\dt v$ in order to mean directed toward~$v$, and $\da v$ to mean directed away from~$v$\sv{.}\lv{; see Figure~\ref{picto}. We will also use the previous definitions with a corner instead of a vertex: in this case, the vertex in question will be the one corresponding to the corner.

\begin{figure}[ht]
		\psfrag{v}[][]{\textcolor{violet}{$v$}}
		\psfrag{h}[][]{$\textcolor{blue}{h}\dt \textcolor{violet}{v}$}
		\psfrag{i}[r][r]{$\textcolor{vert}{h'}\da \textcolor{violet}{v}$}
	\centering\includegraphics[width=8cm]{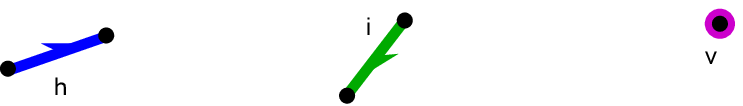}
	\caption{Pictograph representation of a half-edge~$\textcolor{blue}{h}$ directed toward a vertex~$\textcolor{violet}{v}$, and a half-edge~$\textcolor{vert}{h'}$ directed away from~$\textcolor{violet}{v}$.}
	\label{picto}
\end{figure}
}

We end this section by mentioning the following useful elementary facts on bipartite and quasibipartite plane maps. \sv{See the extended version of this paper for a proof.}

\begin{prop}\hspace{0pt}\label{propbq}
The following holds.
\begin{enumerate}[label=(\textit{\roman*})]
	\item In a bipartite map, no edge can be parallel to a vertex. More precisely, for any given face and any given vertex, exactly half of the half-edges incident to the face are directed toward the vertex, the other half being directed away from the vertex.\label{pi}
	\item In a quasibipartite map, a cycle has odd length if and only if it separates the two odd-degree faces\footnote{Recall that, by the Jordan Curve Theorem, a cycle in a plane map always separates the map into exactly two connected components.}. Moreover, for any given vertex~$v$, among the~$a$ half-edges incident to an odd-degree face, exactly one is parallel to~$v$, $(a-1)/2$ are directed toward~$v$ and $(a-1)/2$ are directed away from~$v$.\lv{ For an even-degree face, either zero or two of its incident half-edges are parallel to~$v$, and the remaining incident half-edges are evenly split between those that are directed toward~$v$ and those that are directed away from~$v$.}\label{pii}
\end{enumerate}
\end{prop}

\lv{
\begin{pre}
It comes from the fact that the number of odd-degree faces in any map must be even (recall that the sum of the face degrees is even as it equals twice the number of edges). This implies that, in any plane map, each of the two components separated by a cycle contains an odd number of odd-degree faces if and only if the cycle has odd length (as such a component amounts to a map whose faces are those of the original map that belong to the component plus one face whose degree is the length of the cycle). The first statement of~\ref{pii} follows from this observation, as well as the fact that, in a bipartite map, each cycle has even length.  

Let us now consider an edge~$e$ parallel to a vertex~$v$ in some map. Consider for each extremity of~$e$ a geodesic to the vertex~$v$ and only keep the parts of these two geodesics linking the extremities of~$e$ to their first meeting point. Concatenating these two geodesics together and with~$e$, we obtain some cycle. As the edge is parallel to the vertex, both geodesics have the same length so that the cycle has odd length. From the previous observation, this situation cannot happen in a bipartite map, so that the first statement of~\ref{pi} follows.

As a consequence, in a bipartite map, the distances from some vertex~$v$ to the vertices encountered when traveling along the boundary of a degree $a$ face form an $a$-step bridge whose steps are either~$+1$ or~$-1$. As a result, exactly half of them are~$-1$ steps, that is, half of the incident half-edges are directed toward~$v$.

The same argument as above shows that, in any map, the half-edges incident to a given face that are not parallel to a given vertex are evenly split between those that are directed toward the vertex and those that are directed away from the vertex. Moreover, for obvious parity reasons, the number of half-edges incident to a given face that are parallel to a given vertex has the same parity as the degree of the face. It thus remains to show that, in a quasibipartite map, no more than two edges incident to a given face can be parallel to some given vertex. In fact, consider that we have a map with a face~$f$ having~$k$ edges incident to it that are parallel to some given vertex~$v$. Let us denote by~$e_1$, \ldots, $e_k$ these edges, arranged in counterclockwise order around~$f$. Then consider, for each extremity of each of these~$k$ edges, a geodesic from~$v$. Up to changing the geodesics, one may suppose that any two geodesics never meet again after the time they split. These~$2k$ geodesics are thus arranged into a tree structure with~$2k$ leaves such that, read from left to right, the $2i-1$-th and $2i$-th leaves are the extremities of~$e_i$, for $1\le i \le k$. Using these geodesics, we obtain as above one odd-length cycle per edge~$e_i$, $1\le i \le k$ and the connected components delimited by these cycles that do not contain~$f$ are pairwise disjoint. As a result, each of these~$k$ components must contain an odd number of faces, so that~$k$ is smaller than the number of odd-degree faces of the map, which is~$2$ in a quasibipartite map.
\end{pre}
}

\section{Adding two corners to a face in a bipartite map}\label{secbb}

Throughout this section, we fix an $r$-tuple $\ba=(a_1,\ldots, a_r)$ of positive even integers and we define $\tilde \ba\de( a_1+2, a_2,\ldots, a_r)$ as in the statement of Proposition~\ref{propp2}. We consider on the one hand the set~$\cM$ of plane maps of type~$\ba$ carrying one distinguished edge and two distinguished corners in the first face. On the other hand, we consider the set~$\tilde\cM$ of plane maps of type~$\tilde \ba$ carrying one distinguished vertex and two different distinguished half-edges incident to the first face, and that are both directed toward the distinguished vertex.\lv{ The following pictograph summarizes our definitions (the red $+2$ on the right means that ths size of~$f_1$ has increased by~$2$, and the red arrowhead is the marked corner of~$f_1$):}

\sv{
\begin{center}
		\psfrag{f}[][][.8]{$f_1$}
		\psfrag{2}[][][.7]{\textcolor{red}{$+2$}}
		\psfrag{v}[][][.8]{\textcolor{violet}{$v$}}
		\psfrag{M}[][][.8]{$\cM$}
		\psfrag{N}[][][.8]{$\tilde\cM$}
		\psfrag{c}[][][.8]{\textcolor{blue}{$c$}}
		\psfrag{d}[][][.8]{\textcolor{vert}{$c'$}}
		\psfrag{h}[][][.8]{$\textcolor{blue}{h}\dt \textcolor{violet}{v}$}
		\psfrag{i}[r][r][.8]{$\textcolor{vert}{h'}\dt \textcolor{violet}{v}$}
		\psfrag{e}[][][.8]{\textcolor{violet}{$e$}}
	\centering\includegraphics[width=13cm]{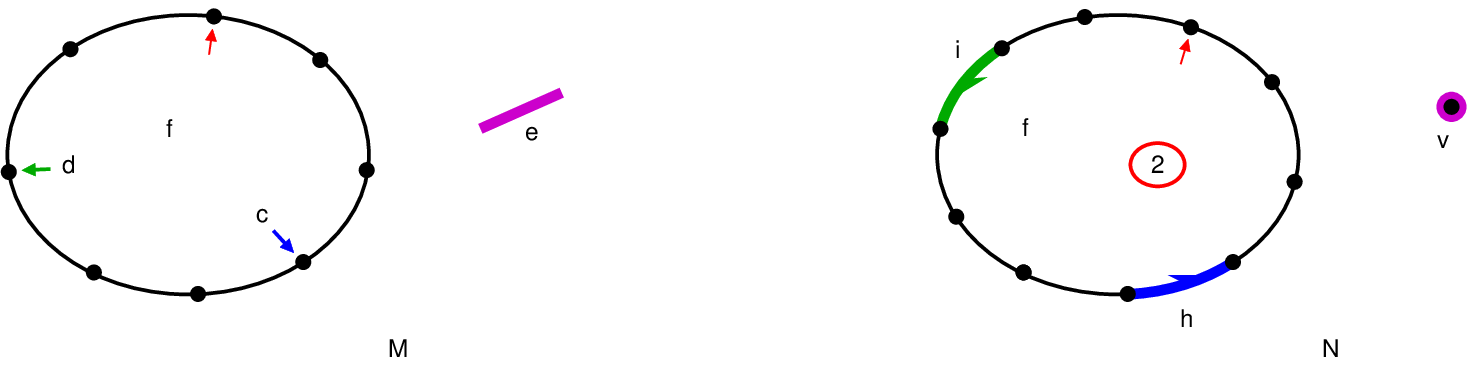}
\end{center}
}\lv{
\begin{center}
		\psfrag{f}[][]{$f_1$}
		\psfrag{2}[][][.9]{\textcolor{red}{$+2$}}
		\psfrag{v}[][][.8]{\textcolor{violet}{$v$}}
		\psfrag{M}[][]{$\cM$}
		\psfrag{N}[][]{$\tilde\cM$}
		\psfrag{c}[][]{\textcolor{blue}{$c$}}
		\psfrag{d}[][]{\textcolor{vert}{$c'$}}
		\psfrag{h}[][]{$\textcolor{blue}{h}\dt \textcolor{violet}{v}$}
		\psfrag{i}[r][r]{$\textcolor{vert}{h'}\dt \textcolor{violet}{v}$}
		\psfrag{e}[][]{\textcolor{violet}{$e$}}
	\centering\includegraphics[width=.95\linewidth]{pictobb}
\end{center}
}

Using Proposition~\ref{propbq}.\ref{pi}, we see that the cardinalities of the sets~$\cM$ and~$\tilde\cM$ are exactly the sides of~\eqref{eqp2}. We now present an explicit bijection between these two sets; this provides a combinatorial interpretation of Proposition~\ref{propp2}. Our bijection is a straightforward generalization of~\cite[Section~4]{bettinelli14ifq}.\lv{ In fact, we treated in the latter reference the case where $a_2=a_3=\ldots=a_r=4$ (up to the irrelevant ordering and corner markings of~$f_2$, \ldots, $f_r$) but the general case can be treated in a similar fashion; everything can be copied almost verbatim. The fact that the faces were of degree~$4$ never intervened; only the fact that the maps were bipartite was of crucial importance. For this reason, we briefly present the construction and refer the interested reader to~\cite[Section~4]{bettinelli14ifq} for more details.

\begin{note}
Although this is not completely obvious, in the case $a_1=2$, $a_2=2p$ and $a_3=a_4=\ldots=a_r=4$, we recover the bijection from \cite[Section~5]{bettinelli14ifq} (still up to face ordering and corner markings). To see this, notice that a map with a distinguished edge corresponds to a map with a distinguished $2$-face by slitting it along the edge. There are then three ways of choosing a $2$-set of corners in this $2$-face; this gives the left hand-side of~\cite[(3)]{bettinelli14ifq}. On the other side, there is only one way of choosing a $2$-set of two different half-edges incident to a $4$-face that are both directed toward a given vertex; we recover the right hand-side of~\cite[(3)]{bettinelli14ifq}. In this setting, it can then be checked that both mappings are indeed the same one.
\end{note}
}

\begin{svc}
\paragraph{Increasing the size.}
\end{svc}
\begin{lvc}
\subsection{Increasing the size}\label{secp2+}
\end{lvc}

Let $(\m;e,c,c')\in \cM$. As~$\m$ is bipartite, $e$ cannot be parallel to~$c$: we denote by~$\vec e $ the corresponding half-edge that is directed toward~$c$, and by~$\pc$ the rightmost geodesic from~$\vec e$ to~$c$. Let us first suppose that $\rev(\vec e)$ is directed toward~$c'$: in this case, the quadruple $(\m;e,c,c')$ is called \emph{simple}. We denote by~$\pc'$ the rightmost geodesic from~$\rev(\vec e)$ to~$c'$ and define the self-avoiding path \sv{$\pp \de \rev(\pc)\ooo \rev(\vec e) \ooo \pc'$. W}\lv{$$\pp \de \rev(\pc)\ooo \rev(\vec e) \ooo \pc'.$$ W}e slit~$\m$ along~$\pp$ from~$c$ to~$c'$, and we denote by~$\pl$ and~$\pr$ the left and right copies of~$\pp$ in the resulting maps. We then sew back $\pl_{1\to [\pp]-1}$ onto $\pr_{2\to [\pp]}$, in the sense that we identify~$\pl_k$ with~$\pr_{k+1}$ for $1\le k \le [\pp]-1$. We denote by~$\tilde\m$ the resulting map and let the outcome of the construction be the quadruple $(\tilde\m;\pl_{[\pc]}^+,\pr_1,\rev(\pl)_{1})$. See Figure~\ref{bb}.

\sv{
\begin{figure}[ht!]
		\psfrag{f}[][][.6]{$f_1$}
		\psfrag{v}[][][.6]{\textcolor{violet}{$v$}}
		\psfrag{a}[][][.6]{\textcolor{red}{$\pl$}}
		\psfrag{b}[][][.6]{\textcolor{red}{$\pr$}}
		\psfrag{m}[][][.6]{$\m$}
		\psfrag{n}[][][.6]{$\tilde\m$}
		\psfrag{c}[][][.6]{\textcolor{blue}{$c$}}
		\psfrag{d}[l][l][.6]{\textcolor{vert}{$c'$}}
		\psfrag{h}[][][.6]{\textcolor{blue}{$h$}}
		\psfrag{i}[][][.6]{\textcolor{vert}{$h'$}}
		\psfrag{e}[][][.6]{\textcolor{violet}{$e$}}
		\psfrag{1}[][][.5]{\textcolor{red}{$\pp_1$}}
		\psfrag{2}[][][.5]{\textcolor{red}{$\pp_2$}}
		\psfrag{3}[][][.5]{\textcolor{red}{$\pp_3$}}
		\psfrag{4}[][][.5]{\textcolor{red}{$\pp_4$}}
		\psfrag{5}[][][.5]{\textcolor{red}{$\pp_5$}}
		\psfrag{6}[][][.5]{\textcolor{red}{$\pp_6$}}
		\psfrag{7}[][][.5]{\textcolor{red}{$\pp_7$}}
		\psfrag{8}[][][.5]{\textcolor{red}{$\pp_8$}}
	\centering\includegraphics[width=.95\linewidth]{bbsv}
	\caption{The mapping from $\cM$ to $\tilde \cM$ in the simple case. We define the path~$\pp$, slit it and sew back after slightly sliding. Only the marked corner of~$f_1$ is represented.}
	\label{bb}
\end{figure}
}\lv{
\begin{figure}[ht!]
		\psfrag{f}[][]{$f_1$}
		\psfrag{v}[][][.8]{\textcolor{violet}{$v$}}
		\psfrag{a}[][][.8]{\textcolor{red}{$\pl$}}
		\psfrag{b}[][][.8]{\textcolor{red}{$\pr$}}
		\psfrag{m}[][]{$\m$}
		\psfrag{n}[][]{$\tilde\m$}
		\psfrag{c}[][]{\textcolor{blue}{$c$}}
		\psfrag{d}[][]{\textcolor{vert}{$c'$}}
		\psfrag{h}[][]{\textcolor{blue}{$h$}}
		\psfrag{i}[][]{\textcolor{vert}{$h'$}}
		\psfrag{e}[][]{\textcolor{violet}{$e$}}
		\psfrag{1}[][][.85]{\textcolor{red}{$\pp_1$}}
		\psfrag{2}[][][.85]{\textcolor{red}{$\pp_2$}}
		\psfrag{3}[][][.85]{\textcolor{red}{$\pp_3$}}
		\psfrag{4}[][][.85]{\textcolor{red}{$\pp_4$}}
		\psfrag{5}[][][.85]{\textcolor{red}{$\pp_5$}}
		\psfrag{6}[][][.85]{\textcolor{red}{$\pp_6$}}
		\psfrag{7}[][][.85]{\textcolor{red}{$\pp_7$}}
		\psfrag{8}[][][.85]{\textcolor{red}{$\pp_8$}}
	\centering\includegraphics[width=11cm]{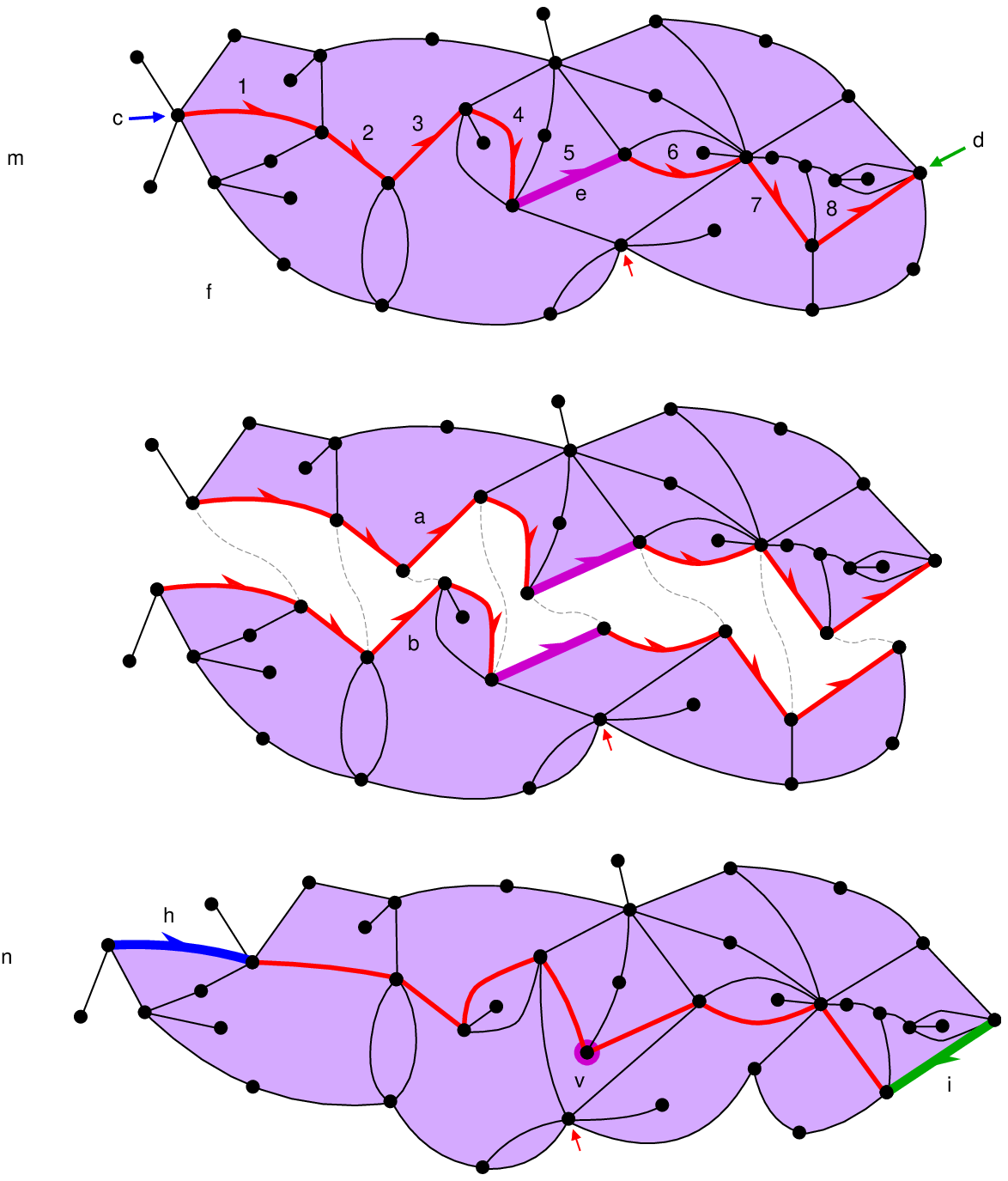}
	\caption{The mapping from $\cM$ to $\tilde \cM$ in the simple case. We define the path~$\pp$, slit it and sew back after slightly sliding. On this picture, only the marked corner of~$f_1$ is represented.}
	\label{bb}
\end{figure}
}

Let us now treat the case where~$\vec e$ is directed toward~$c'$. We denote by~$\pc'$ the rightmost geodesic from~$\vec e$ to~$c'$ and by $i\ge 1$ the smallest integer such that $\pc_i \neq \pc'_i$. As~$\pc$ and~$\pc'$ are rightmost geodesics, we must have $\{\pc_i^+,\dots,\pc_{[\pc]}^+\}\cap\{\pc_i'^+,\dots, \pc_{[\pc']}'^+\}=\varnothing$. The path
$$\pp \de \rev(\pc) \ooo \rev(\vec e) \ooo \vec e \ooo \pc'$$
is thus composed of the self-avoiding path $\rev(\pc_{i\to [\pc]}) \ooo \pc'_{i \to [\pc']}$ together with the self-avoiding path $\vec e \ooo \pc_{1\to i-1}$ (visited twice, first backwards then forward), grafted either to its left or to its right. We say that the path~$\pp$ and the quadruple $(\m;e,c,c')$ are \emph{left-pinched} or \emph{right-pinched} accordingly.

As above, we slit~$\m$ along~$\pp$ from~$c$ to~$c'$, circumventing the pinched part. This still splits~$\m$ into two submaps with a copy of~$\pp$ on the boundary of each but, this time, one copy is a self-avoiding path while the other copy goes back and forth along a ``dangling'' chain of~$i$ edges at some point. We still denote the left and right copies of~$\pp$ by~$\pl$ and~$\pr$ and sew back $\pl_{1\to [\pp]-1}$ onto $\pr_{2\to [\pp]}$. We denote by~$\tilde\m$ the resulting map and let the outcome of the construction be the quadruple $(\tilde\m;\pl_{[\pc]}^+,\pr_1,\rev(\pl)_{1})$ in the left-pinched case and $(\tilde\m;(\rev(\pr))_{[\pc']}^+,\pr_1,\rev(\pl)_{1})$ in the right-pinched case (so that the distinguished vertex is always the tip of the dangling chain).\lv{ See Figure~\ref{bbpinch}. 

\begin{figure}[ht]
		\psfrag{f}[][][.8]{$f_1$}
		\psfrag{v}[][][.8]{\textcolor{violet}{$v$}}
		\psfrag{a}[][][.8]{\textcolor{red}{$\pl$}}
		\psfrag{b}[][][.8]{\textcolor{red}{$\pr$}}
		\psfrag{m}[][][.8]{$\m$}
		\psfrag{n}[][][.8]{$\tilde\m$}
		\psfrag{c}[][]{\textcolor{blue}{$c$}}
		\psfrag{d}[][]{\textcolor{vert}{$c'$}}
		\psfrag{h}[][][.8]{\textcolor{blue}{$h$}}
		\psfrag{i}[][][.8]{\textcolor{vert}{$h'$}}
		\psfrag{e}[][][.8]{\textcolor{violet}{$e$}}		
	\centering\includegraphics[width=.95\linewidth]{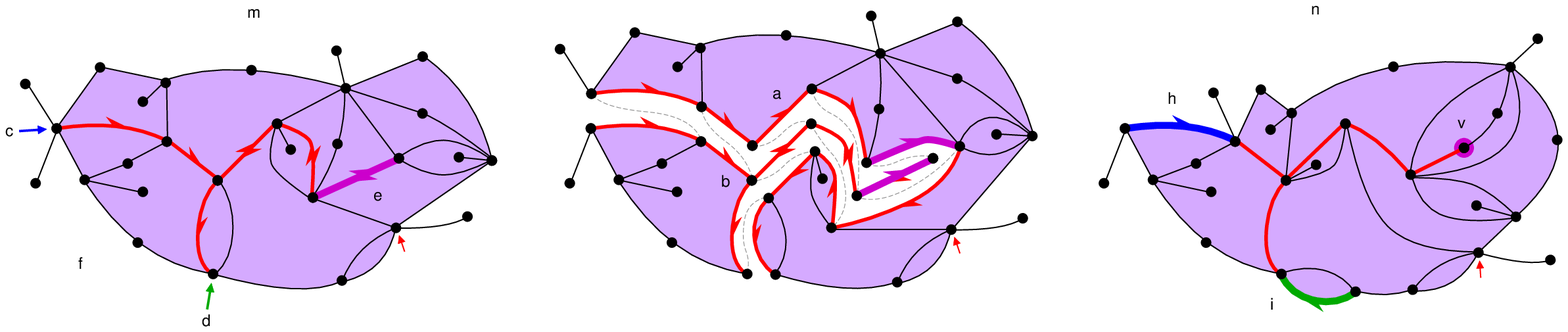}
	\caption{The mapping from $\cM$ to $\tilde \cM$ in the pinched case. We slide along the path~$\pp$ and circumvent its pinched part.}
	\label{bbpinch}
\end{figure}
}

\begin{svc}
\paragraph{Decreasing the size.}
\end{svc}
\begin{lvc}
\subsection{Decreasing the size}\label{secp2-}
\end{lvc}

The inverse mapping takes a quadruple $(\tilde\m;v,h,h') \in \tilde \cM$ and goes as follows. We consider the corner~$h_0$ delimited by~$h$ and its predecessor in the contour of the first face of~$\tilde\m$, and denote by~$\ph$ the leftmost geodesic from this corner to~$v$. As~$h$ is directed toward~$v$, we have that $[\ph]\ge 1$ and $\ph_1=h$. We define~$h_0'$ and~$\ph'$ in a similar fashion with~$h'$ instead of~$h$. Depending on whether~$\ph$ and~$\ph'$ meet before reaching~$v$ or not, the path \sv{$\pp' \de \ph \ooo \rev(\ph')$ i}\lv{$$\pp' \de \ph \ooo \rev(\ph')$$ i}s either self-avoiding or pinched in the sense of the previous \sv{paragraph}\lv{section}.\lv{ To see this, observe that, if~$\ph$ and~$\ph'$ meet before reaching~$v$, denoting by~$i$ and~$j$ the smallest integers such that $\ph_{i}^+=\ph_{j}'^+$, the path $\ph_{1\to i}\ooo \rev(\ph'_{1\to j})$ separates the map into two disjoint components and~$v$ belongs to only one of them.} The quadruple $(\tilde\m;v,h,h')$ is called \emph{simple}, \emph{left-pinched} or \emph{right-pinched} accordingly. We slit~$\tilde\m$ along~$\pp'$ from~$h_0$ to~$h_0'$, denote by~$\pl'$ and~$\pr'$ the left and right copies of~$\pp'$ in the resulting maps and sew $\pl'_{2\to [\pp']}$ onto $\pr'_{1\to [\pp']-1}$. In the resulting map, $\pl'_1$ and $(\rev(\pr'))_1$ are dangling edges. We suppress them and denote respectively by~$c$ and~$c'$ the corners they define. We denote by~$\m$ the map we finally obtain and let the outcome of the construction be the quadruple $(\m;e,c,c')$, where~$e$ is the edge corresponding to~$\pl'_{[\ph]+1}$.

\begin{svc}
\paragraph{The previous mappings are inverse one from another.}
\end{svc}
\begin{lvc}
\subsection{The previous mappings are inverse one from another}\label{secbbinv}
\end{lvc}

In fact, through the mappings of the two previous \sv{paragraphs}\lv{sections}, simple quadruples correspond to simple quadruples, left-pinched quadruples correspond to left-pinched quadruples and right-pinched quadruples correspond to right-pinched quadruples.

The proof that the previous mappings are inverse one from another can be copied almost verbatim from \cite[Proof of Theorem~3]{bettinelli14ifq}.\lv{ For the sake of self-containment, we will very briefly recall the main steps of this proof.} Alternatively, we will see in Section~\ref{secdecomp} that these mappings can be seen as compositions of simpler slit-slide-sew bijections; this will provide an alternate, arguably simpler, proof.

\lv{
In the notation of Section~\ref{secp2+}, the map~$\tilde\m$ is clearly of the desired type~$\tilde \ba$. We claim that the image in~$\tilde\m$ of the path $\pr_{1\to [\pc]+1}$ is the leftmost geodesic from the last corner before~$\pr_1$ in the contour of the first face toward~$\pr_{[\pc]+1}^+=\pl_{[\pc]}^+$. As the construction is completely symmetric, this will also entail that $\rev(\pl)_{1\to [\pc']+1}$ is the leftmost geodesic from the last corner before~$\rev(\pl)_1$ to~$\rev(\pl)_{[\pc']+1}^+$ and, as a result, that $\pp'=\pr \ooo\pl_{[\pp]}$. From this claim, we thus conclude that the outcome of the construction belongs to~$\tilde\cM$ and that the application of the construction of Section~\ref{secp2-} to it gives back the initial quadruple $(\m;e,c,c')$. In order to show this claim, we track back the considered geodesics into the original map and see how they behave. A similar argument in a simpler case will be used during Section~\ref{sectrans}; see Figure~\ref{pftransfer}.
}

\section{Adding one corner to two faces in a bipartite map}\label{secbq}

We now\lv{ use the setting of Proposition~\ref{propp1p1}. Namely, we} fix an $r$-tuple $\ba=(a_1,\ldots, a_r)$ of positive even integers and we define $\tilde \ba\de( a_1+1, a_2+1,a_3,\ldots, a_r)$. We let~$\cM$ be the set of plane maps of type~$\ba$ carrying one distinguished edge, one distinguished corner in the first face and one distinguished corner in the second face. We let~$\tilde\cM$ be the set of plane maps of type~$\tilde \ba$ carrying one distinguished vertex and two distinguished half-edges directed toward it, one being incident to the first face and one being incident to the second face.

\sv{
\begin{center}
		\psfrag{f}[][][.8]{$f_1$}
		\psfrag{g}[][][.8]{$f_2$}
		\psfrag{1}[][][.7]{\textcolor{red}{$+1$}}
		\psfrag{v}[][][.8]{\textcolor{violet}{$v$}}
		\psfrag{M}[][][.8]{$\cM$}
		\psfrag{N}[][][.8]{$\tilde\cM$}
		\psfrag{c}[][][.8]{\textcolor{blue}{$c$}}
		\psfrag{d}[][][.8]{\textcolor{vert}{$c'$}}
		\psfrag{h}[][][.8]{$\textcolor{blue}{h}\dt \textcolor{violet}{v}$}
		\psfrag{i}[][][.8]{$\textcolor{vert}{h'}\dt \textcolor{violet}{v}$}
		\psfrag{e}[][][.8]{\textcolor{violet}{$e$}}
	\centering\includegraphics[width=13cm]{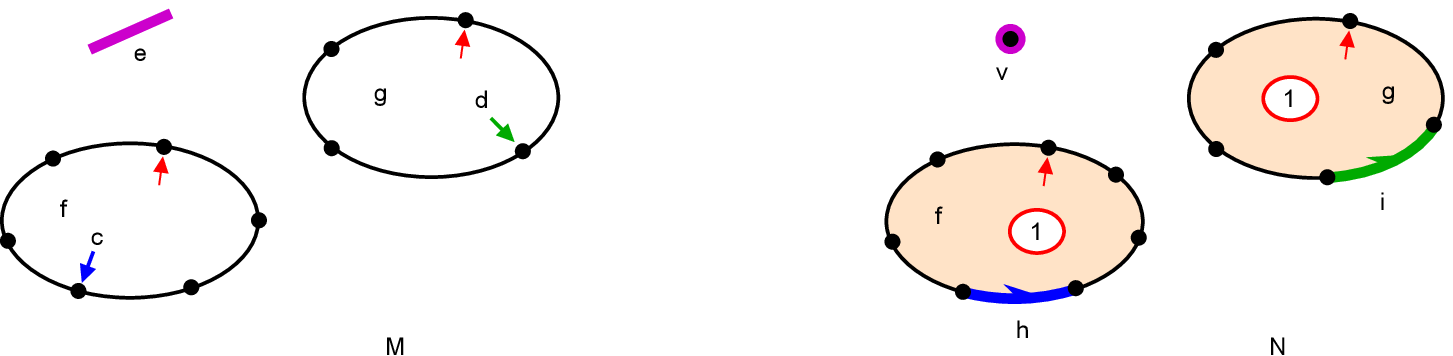}
\end{center}
}\lv{
\begin{center}
		\psfrag{f}[][]{$f_1$}
		\psfrag{g}[][]{$f_2$}
		\psfrag{1}[][][.9]{\textcolor{red}{$+1$}}
		\psfrag{v}[][][.8]{\textcolor{violet}{$v$}}
		\psfrag{M}[][]{$\cM$}
		\psfrag{N}[][]{$\tilde\cM$}
		\psfrag{c}[][]{\textcolor{blue}{$c$}}
		\psfrag{d}[][]{\textcolor{vert}{$c'$}}
		\psfrag{h}[][]{$\textcolor{blue}{h}\dt \textcolor{violet}{v}$}
		\psfrag{i}[][]{$\textcolor{vert}{h'}\dt \textcolor{violet}{v}$}
		\psfrag{e}[][]{\textcolor{violet}{$e$}}
	\centering\includegraphics[width=.95\linewidth]{pictobq}
\end{center}
}

\begin{figure}[ht]
		\psfrag{f}[][][.6]{$f_1$}
		\psfrag{g}[][][.6]{$f_2$}
		\psfrag{v}[][][.6]{\textcolor{violet}{$v$}}
		\psfrag{a}[][][.6]{\textcolor{red}{$\pl$}}
		\psfrag{b}[][][.6]{\textcolor{red}{$\pr$}}
		\psfrag{m}[][][.6]{$\m$}
		\psfrag{n}[][][.6]{$\tilde\m$}
		\psfrag{c}[][][.6]{\textcolor{blue}{$c$}}
		\psfrag{d}[][][.6]{\textcolor{vert}{$c'$}}
		\psfrag{h}[][][.6]{\textcolor{blue}{$h$}}
		\psfrag{i}[][][.6]{\textcolor{vert}{$h'$}}
		\psfrag{e}[][][.6]{\textcolor{violet}{$e$}}
	\centering\includegraphics[width=.95\linewidth]{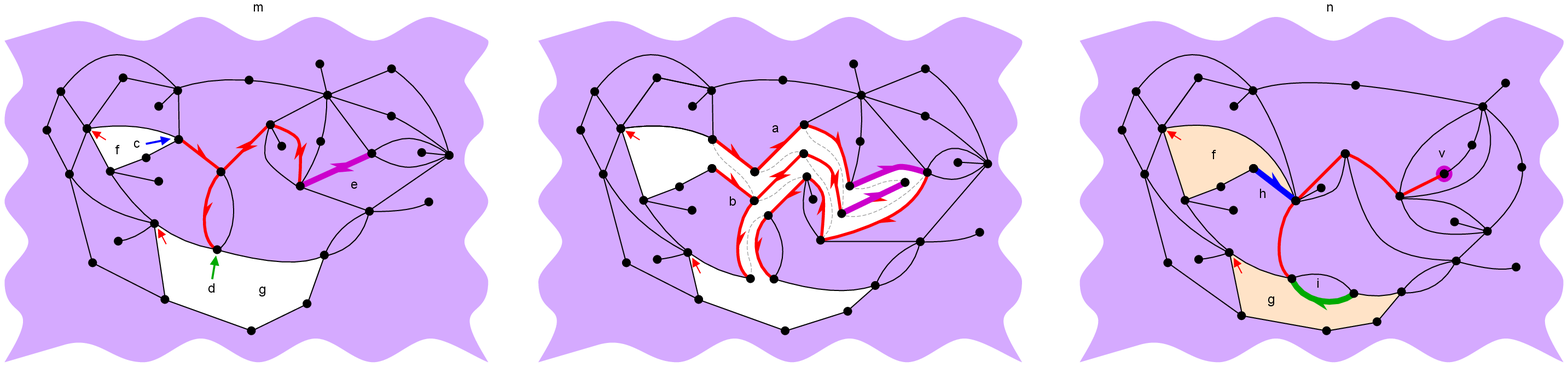}
	\caption{The mapping from $\cM$ to $\tilde \cM$ in the pinched case.}
	\label{bqpinch}
\end{figure}

The cardinality of~$\cM$ is clearly equal to the left-hand side of~\eqref{eqp1p1} and we see that the cardinality of~$\tilde\cM$ is equal to the right-hand side of~\eqref{eqp1p1} by using Proposition~\ref{propbq}.\ref{pii}. The mappings interpreting Proposition~\ref{propp1p1} are described exactly as in \sv{the previous section}\lv{Sections~\ref{secp2+} and~\ref{secp2-}}: see Figure~\ref{bqpinch}. The only difference is that the paths~$\pp$ and~$\pp'$ no longer disconnect the maps\lv{ (recall Remark~\ref{remdiscon})}; this bears no effects in the description of the mappings.

\lv{In order to see that these mappings are well defined, one only needs}\sv{It is not very hard} to see that the paths~$\pp$ and~$\pp'$ are as before (self-avoiding or pinched)\sv{; we refer the reader to the extended version.}\lv{. For~$\pp$, the arguments we used in Section~\ref{secbb} still hold. For~$\pp'$, one needs an extra argument when~$\ph$ and~$\ph'$ meet before reaching~$v$. Let~$i$ and~$i'$ be the smallest integers such that $\ph_{i}^+=\ph_{i'}'^+$ and assume by contradiction that $\ph_{i+1}\neq \ph_{i'+1}'$. As~$\ph$ and~$\ph'$ are leftmost geodesics toward~$v$, they are bound to meet again and the half-edges~$\ph_i$, $\ph_{i+1}$, $\ph_{i'}'$, $\ph'_{i'+1}$ are arranged in counterclockwise order around~$\ph_{i}^+$. Let~$j$ and~$j'$ be the second smallest integers such that $\ph_{j}^+=\ph_{j'}'^+$. As~$\ph$ and~$\ph'$ are geodesics, we have that $j-i=j'-i'$. Now, the path $\ph_{i+1\to j}\ooo \rev(\ph'_{i'+1\to j'})$ is an even-length cycle that cannot intersect $\ph_{1\to i}\ooo \rev(\ph'_{1\to i'})$ except at~$\ph_{i}^+$. As a result, this cycle separates~$f_1$ from~$f_2$, a contradiction to Proposition~\ref{propbq}.\ref{pii}.} 

\lv{
The proof that these mappings are inverse one from another goes almost exactly as in Section~\ref{secbbinv}. The only difference is that the maps of~$\tilde\cM$ are no longer bipartite: there might thus be odd-length cycles. Hopefully, thanks to Proposition~\ref{propbq}.\ref{pii}, such cycles do not alter the arguments of the proof, as they have to seperate the odd-degree faces. We do not linger on this technical issue; once again, the decomposition of Section~\ref{sectrans} will provide a simpler proof.
}

\section{Transfer bijections}\label{sectrans}

Let us now see how the previous mappings can be decomposed as two more elementary mappings.\lv{ More precisely, we used an edge~$e$ and two corners~$c$ and~$c'$ (either of the same face or of two different faces) and built a path linking the corners as the concatenation of two paths going from the edge to each of the corners. We will see this operation as the result of two slit-slide-sew bijections as follows. First, we replace the distinguished edge~$e$ with an $r+1$-th face~$f_{r+1}$ of degree~$2$ by doubling the edge (the marked corner of this face is arbitrarily chosen). Next, we subsequently apply two mappings that each transfers a corner from~$f_{r+1}$ to the faces containing~$c$ and~$c'$. As a result, $f_{r+1}$ completely vanishes into a vertex. Let us see in more details these transfer mappings and come back to this decomposition more precisely later on (in Section~\ref{secdecomp}).}

\subsection{Transferring from a face of degree at least two}

We start with \sv{the setting of }Proposition~\ref{propp1m1}.\lv{ Let $\ba=(a_1,\ldots, a_{r+1})$ be an $r+1$-tuple of positive integers with~$a_{r+1}\ge 2$ and
\begin{itemize}
	\item either all even,
	\item or such that only~$a_{r+1}$ and one other coordinate are odd,
\end{itemize}
and let $\tilde \ba=(\tilde a_1,\ldots,\tilde a_{r+1})\de( a_1+1, a_2,\ldots, a_{r}, a_{r+1}-1)$. We interpret the sides of~\eqref{eqp1m1} as follows.} We let~$\cM$ be the set of plane maps of type~$\ba$ carrying one distinguished corner~$c$ in the first face and one distinguished half-edge~$h'$ incident to the $r+1$-th face and directed toward~$c$. We define~$\tilde\cM$ as the set of plane maps of type~$\tilde\ba$ carrying one distinguished corner~$c'$ in~$f_{r+1}$ and one distinguished half-edge~$h$ incident to the first face and directed \emph{away from}~$c'$.

\sv{
\begin{center}
		\psfrag{f}[][][.8]{$f_1$}
		\psfrag{g}[][][.8]{$f_{r+1}$}
		\psfrag{1}[][][.7]{\textcolor{red}{$+1$}}
		\psfrag{9}[][][.7]{\textcolor{red}{$-1$}}
		\psfrag{M}[][][.8]{$\cM$}
		\psfrag{N}[][][.8]{$\tilde\cM$}
		\psfrag{c}[][][.8]{\textcolor{blue}{$c$}}
		\psfrag{d}[][][.8]{\textcolor{vert}{$c'$}}
		\psfrag{h}[][][.8]{$\textcolor{blue}{h}\da \textcolor{vert}{c'}$}
		\psfrag{i}[][][.8]{$\textcolor{vert}{h'}\dt \textcolor{blue}{c}$}
	\centering\includegraphics[width=13cm]{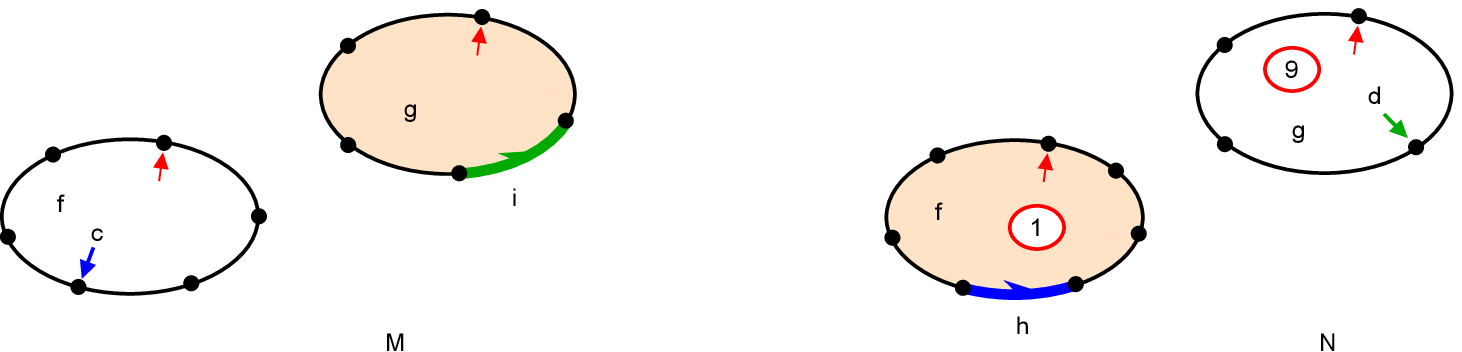}
\end{center}
}\lv{
\begin{center}
		\psfrag{f}[][]{$f_1$}
		\psfrag{g}[][]{$f_{r+1}$}
		\psfrag{1}[][][.9]{\textcolor{red}{$+1$}}
		\psfrag{9}[][][.9]{\textcolor{red}{$-1$}}
		\psfrag{M}[][]{$\cM$}
		\psfrag{N}[][]{$\tilde\cM$}
		\psfrag{c}[][]{\textcolor{blue}{$c$}}
		\psfrag{d}[][]{\textcolor{vert}{$c'$}}
		\psfrag{h}[][]{$\textcolor{blue}{h}\da \textcolor{vert}{c'}$}
		\psfrag{i}[][]{$\textcolor{vert}{h'}\dt \textcolor{blue}{c}$}
	\centering\includegraphics[width=.95\linewidth]{pictoelem}
	
	\vspace{3mm}
	\emph{On the above pictograph, we chose to depict the case of an even~$a_1$ and an odd~$a_{r+1}$.}
\end{center}

Every map of type~$\ba$ has $\big\lfloor {a_{r+1}}/{2}\big\rfloor$ half-edges incident to the $r+1$-th face that are directed toward a given corner. This comes from Proposition~\ref{propbq}.\ref{pi} when every coordinate of~$\ba$ is even, and from Proposition~\ref{propbq}.\ref{pii} when the map is quasibipartite, as~$a_{r+1}$ is odd. As a result, the cardinality of~$\cM$ is the left-hand side of~\eqref{eqp1m1}. The conditions on~$\ba$ imply that either all coordinates of~$\tilde\ba$ are even, or only~$\tilde a_{1}$ and one other coordinate are odd. By the above argument, the cardinality of~$\tilde\cM$ is the right-hand side of~\eqref{eqp1m1}.
}

Let us describe the mappings (see \sv{the left part of }Figure~\ref{transfer2}) between~$\cM$ and~$\tilde \cM$. Let $(\m;c,h')\in \cM$. We consider the corner~$h'_0$ delimited by~$h'$ and its predecessor in the contour of~$f_{r+1}$, and denote by~$\ph'$ the leftmost geodesic from~$h'_0$ to~$c$. We slit~$\m$ along~$\ph'$ from~$h'_0$ to~$c$, denote by~$\pl'$ and~$\pr'$ the left and right copies of~$\ph'$ in the resulting map and sew $\pl'_{2\to [\ph']}$ onto $\pr'_{1\to [\ph']-1}$. In the resulting map, we denote by~$h$ the half-edge~$\pr'_{[\ph']}$, suppress the dangling edge~$\pl'_1$ and denote by~$c'$ the corner it defines. We then denote by~$\tilde\m$ the resulting map and let the outcome of the construction be $\Phi_{\operatorname{left}}(\m;c,h')\de(\tilde\m;c',h)$.

Conversely, starting from $(\tilde\m;c',h)\in \tilde\cM$, we consider the corner~$h_0$ delimited by~$h$ and its successor in the contour of~$f_1$, and denote by~$\ph$ the rightmost geodesic from~$h_0$ to~$c'$. We slit~$\tilde\m$ along~$\ph$ from~$h_0$ to~$c'$, denote by~$\pl$ and~$\pr$ the left and right copies of~$\ph$ in the resulting map and sew $\pl_{1\to [\ph]-1}$ onto $\pr_{2\to [\ph]}$. In the resulting map, we denote by~$h'$ the half-edge~$\rev(\pl)_{1}$, suppress the dangling edge~$\pr_1$ and denote by~$c$ the corner it defines. We then denote by~$\m$ the resulting map and let the outcome of the construction be $\Phi_{\operatorname{right}}(\tilde\m;c',h)\de(\m;c,h')$.

\sv{
\begin{figure}[ht!]
		\psfrag{g}[][][.6]{$f_1$}
		\psfrag{f}[][][.6]{$f_{r+1}$}		
		\psfrag{p}[][][.8]{\textcolor{red}{$\ph'$}}
		\psfrag{l}[][][.8]{\textcolor{red}{$\pl'$}}
		\psfrag{r}[][][.8]{\textcolor{red}{$\pr'$}}
		\psfrag{m}[][][.8]{$\m$}
		\psfrag{n}[][][.8]{$\tilde\m$}
		\psfrag{c}[][][.8]{\textcolor{blue}{$c$}}
		\psfrag{d}[][][.8]{\textcolor{vert}{$c'$}}
		\psfrag{h}[][][.8]{\textcolor{blue}{$h$}}
		\psfrag{i}[][][.8]{\textcolor{vert}{$h'$}}
		\psfrag{e}[][][.8]{\textcolor{vert}{$h'_0$}}
	\centering\includegraphics[height=12cm]{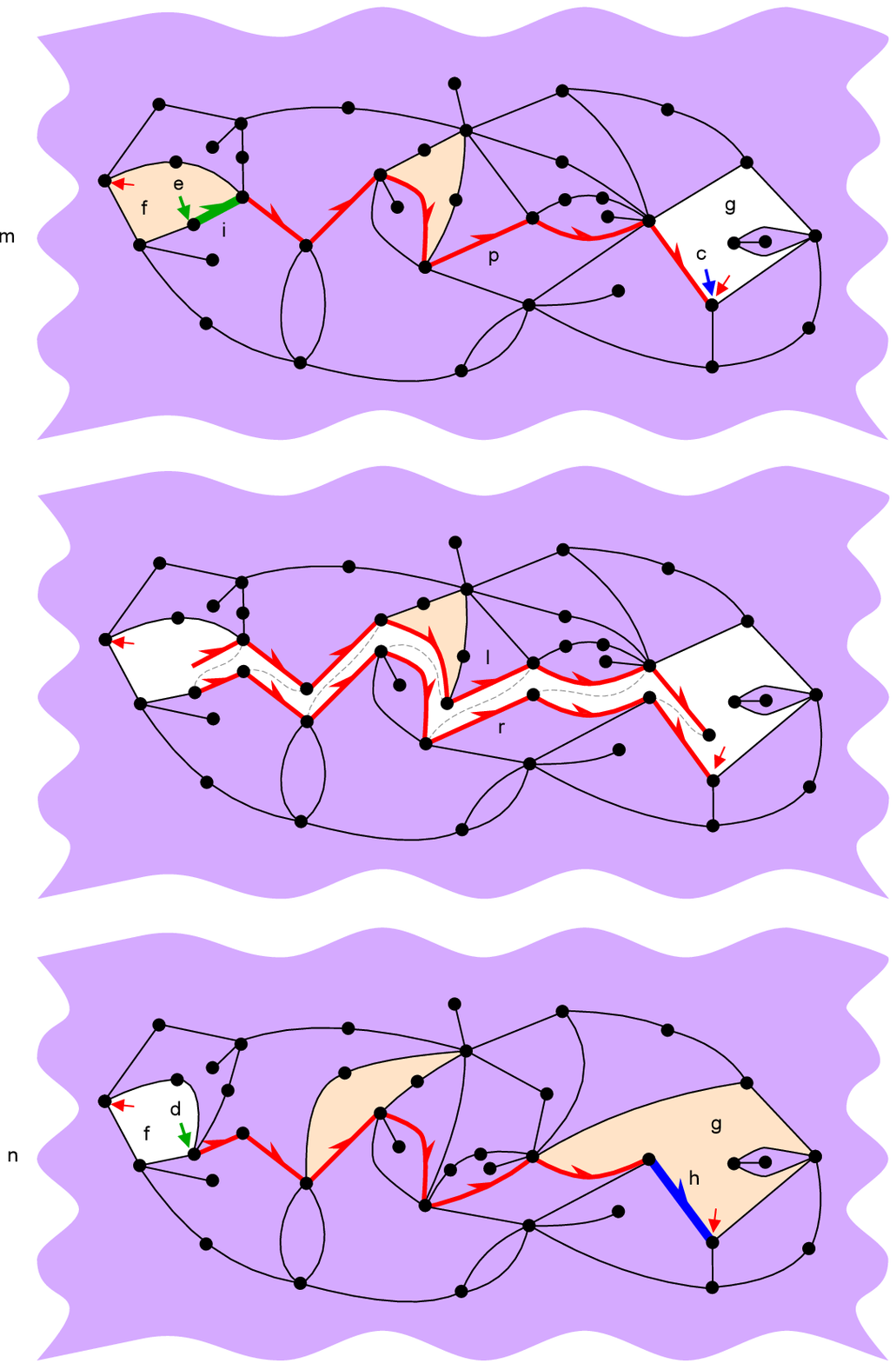}\qquad
		\psfrag{f}[][][.6]{$f_1$}
		\psfrag{g}[][][.6]{$f_{r+1}$}		
		\psfrag{p}[][][.8]{\textcolor{red}{$\pp'$}}
		\psfrag{l}[][][.8]{\textcolor{red}{$\pl'$}}
		\psfrag{r}[][][.8]{\textcolor{red}{$\pr'$}}
		\psfrag{m}[][][.8]{$\m$}
		\psfrag{n}[][][.8]{$\tilde\m$}
		\psfrag{c}[][][.8]{\textcolor{blue}{$c$}}
		\psfrag{h}[][][.8]{\textcolor{blue}{$h$}}
		\psfrag{e}[][][.8]{\textcolor{blue}{$h_0$}}
		\psfrag{v}[][][.8]{\textcolor{violet}{$v$}}
	\centering\includegraphics[height=12cm]{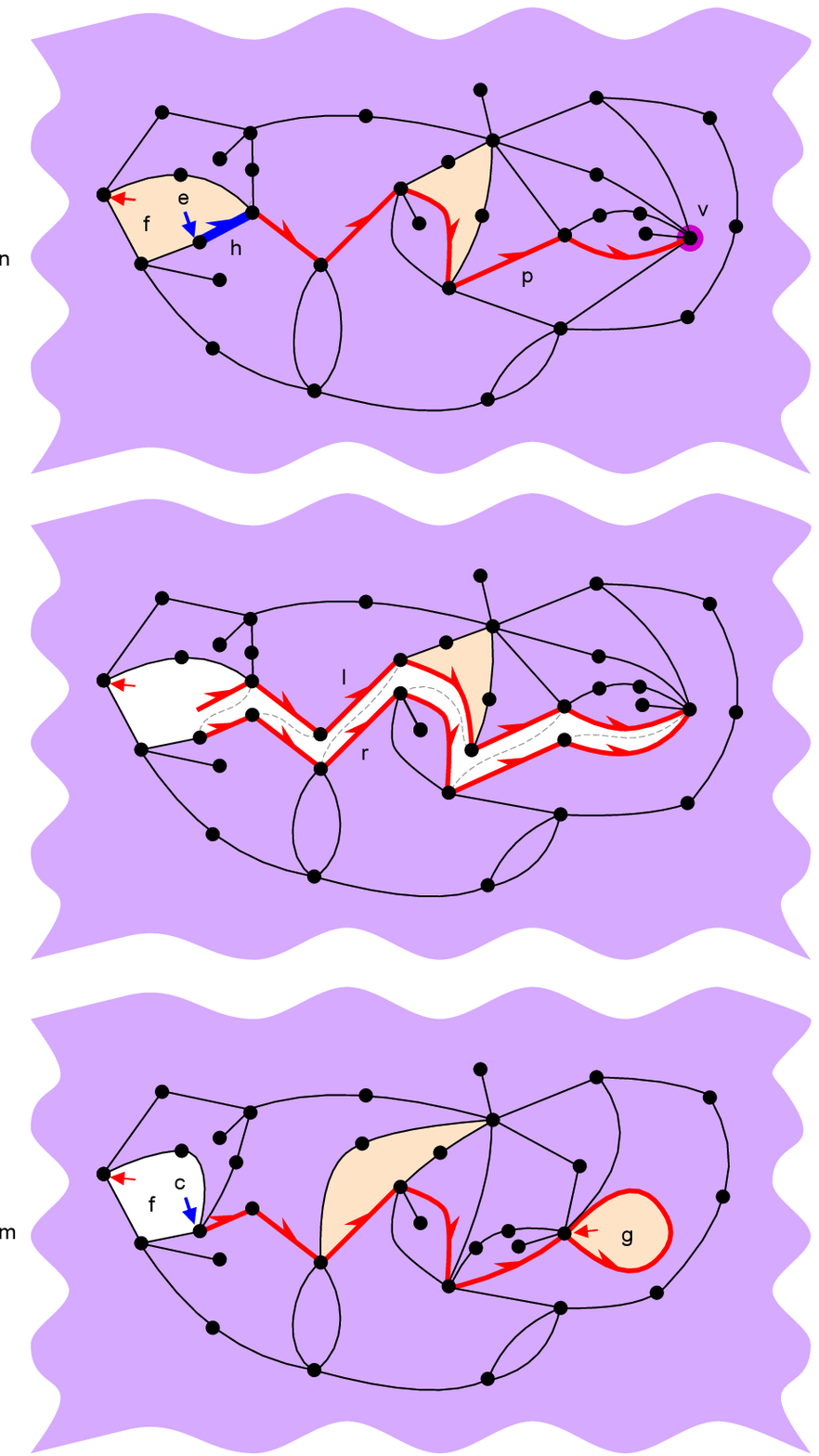}
	\caption{The transfer mappings from a face of degree more than~$2$ (left) or~$1$ (right).}
	\label{transfer2}
\end{figure}
}\lv{
\begin{figure}[ht!]
		\psfrag{g}[][][.8]{$f_1$}
		\psfrag{f}[][][.8]{$f_{r+1}$}		
		\psfrag{p}[][][.8]{\textcolor{red}{$\ph'$}}
		\psfrag{l}[][][.8]{\textcolor{red}{$\pl'$}}
		\psfrag{r}[][][.8]{\textcolor{red}{$\pr'$}}
		\psfrag{m}[][]{$\m$}
		\psfrag{n}[][]{$\tilde\m$}
		\psfrag{c}[][]{\textcolor{blue}{$c$}}
		\psfrag{d}[][]{\textcolor{vert}{$c'$}}
		\psfrag{h}[][]{\textcolor{blue}{$h$}}
		\psfrag{i}[][]{\textcolor{vert}{$h'$}}
		\psfrag{e}[][][.8]{\textcolor{vert}{$h'_0$}}
	\centering\includegraphics[width=10cm]{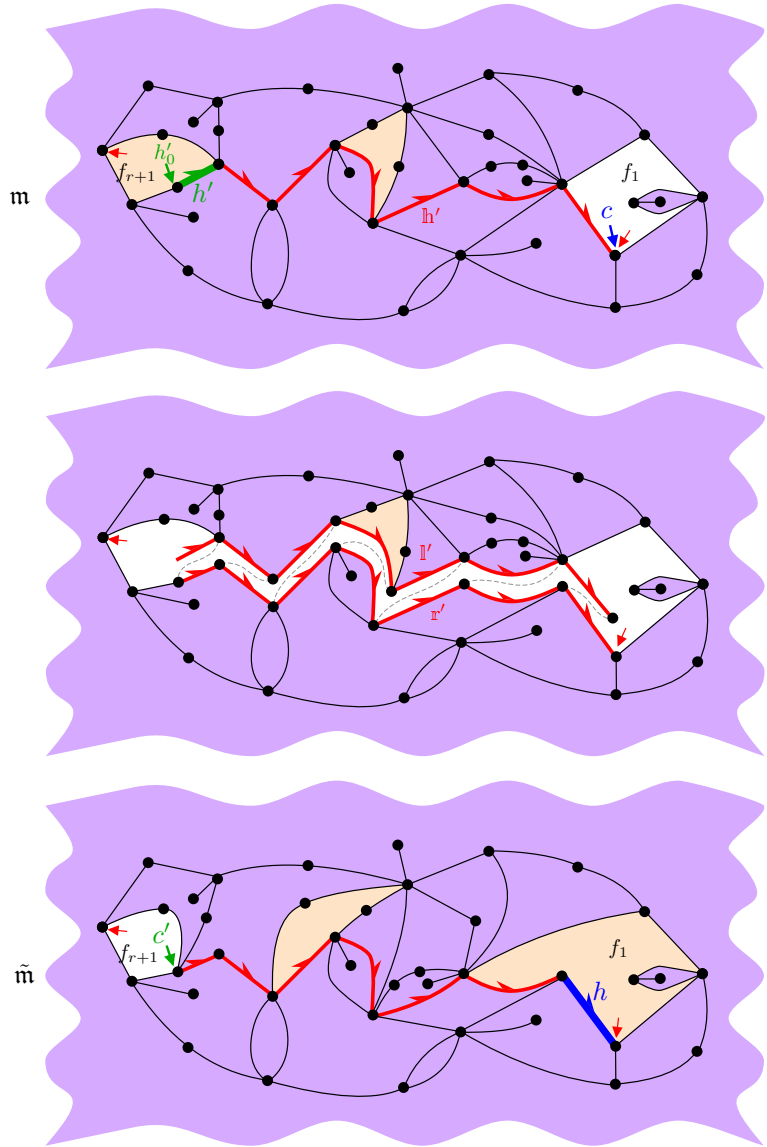}
	\caption{The transfer mappings.}
	\label{transfer2}
\end{figure}
}

\begin{thm}\label{thmtransfer2}
The mappings $\Phi_{\operatorname{left}}:\cM\to\tilde\cM$ and $\Phi_{\operatorname{right}}:\tilde\cM\to\cM$ are inverse bijections.
\end{thm}

\sv{We refer the reader to the extended version for the proof.}
\lv{
\begin{pre}
It is clear from the constructions that the types of the maps are as desired and that the mappings are inverse one from another, provided the conditions on the distinguished half-edges are satisfied and the sliding paths correspond. 

We consider $(\m;c,h')\in \cM$ and define $(\tilde\m;c',h)\de\Phi_{\operatorname{left}}(\m;c,h')$. Let us see that the image in~$\tilde\m$ of the path~$\rev(\pr')$ is the path~$\ph$. This will entail that~$h$ is directed away from~$c'$, so that $(\tilde\m;c',h)\in\tilde\cM$, and that $\Phi_{\operatorname{right}}\circ\Phi_{\operatorname{left}}$ is the identity on~$\cM$.

We argue by contradiction and suppose that $\ph\neq\rev(\pr')$ (we keep the notation~$\rev(\pr')$ for its image in~$\tilde\m$). Then~$\ph$ has to leave the path~$\rev(\pr')$ at some point (to its left or to its right) and come back to it at some other point (from its left or from its right). It is easy to check that these four possibilities contradict the fact that~$\ph'$ is the leftmost geodesic from~$h'_0$ to~$c$\,; see Figure~\ref{pftransfer}.

\begin{figure}[ht]
		\psfrag{r}[][][.8]{\textcolor{red}{$\pr'$}}
		\psfrag{t}[][][.8]{\textcolor{red}{$\ph'$}}
		
		\psfrag{m}[][]{$\m$}
		\psfrag{n}[][]{$\tilde\m$}
		\psfrag{c}[][][.8]{\textcolor{blue}{$c$}}
		\psfrag{d}[][][.8]{\textcolor{vert}{$c'$}}
		\psfrag{h}[][][.8]{\textcolor{blue}{$h$}}
		\psfrag{i}[][][.8]{\textcolor{vert}{$h'$}}
		\psfrag{j}[][][.8]{\textcolor{vert}{$h'_0$}}
		\psfrag{e}[][][.8]{\textcolor{blue}{$h_0$}}
		
		\psfrag{a}[][][.6]{\textcolor{violet}{$\bm{\le}$}}
		\psfrag{b}[][][.6]{\textcolor{violet}{$\bm{<}$}}
	\centering\includegraphics[width=.95\linewidth]{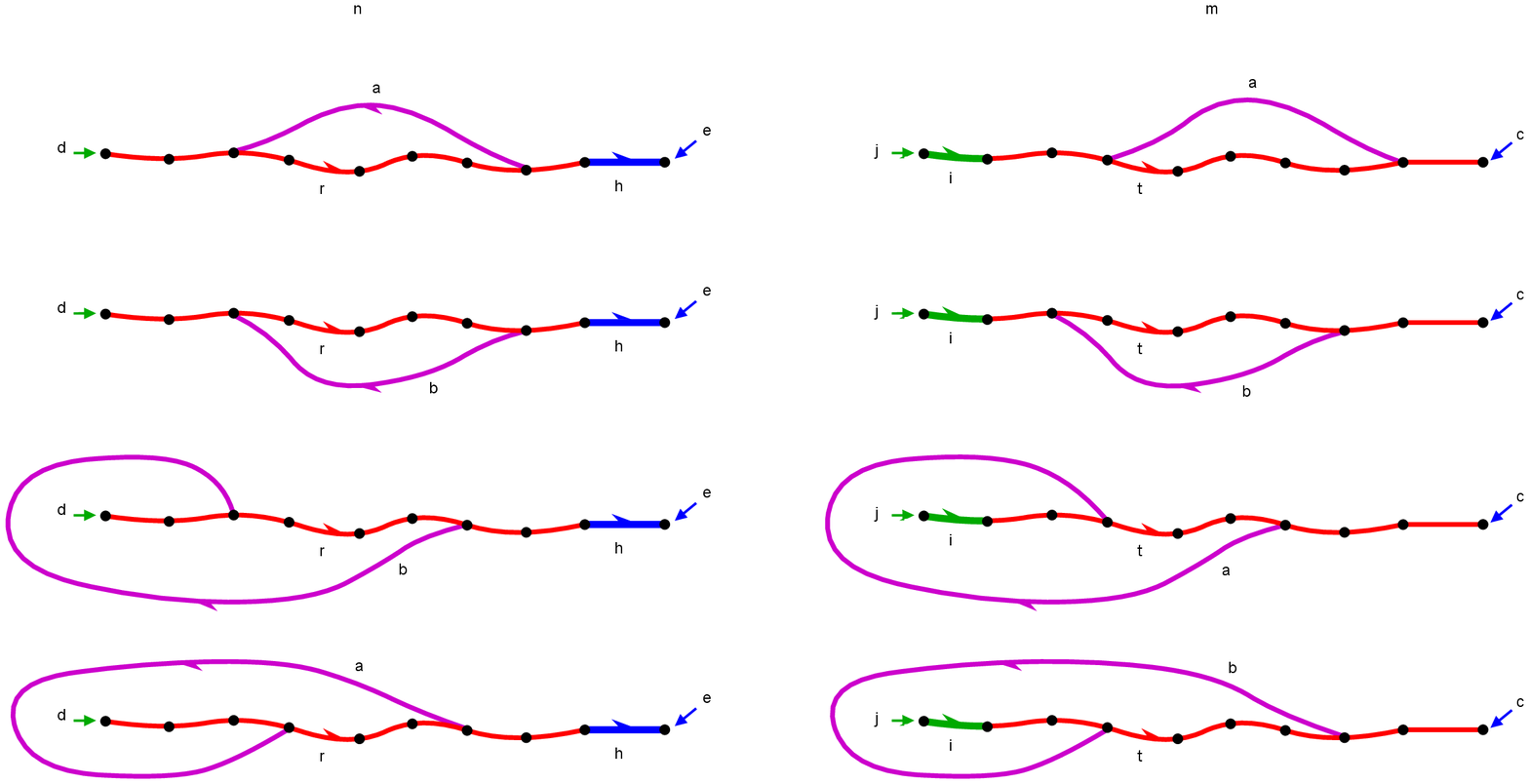}
	\caption{Proof of the fact that the leftmost geodesic~$\ph'$ from~$h'_0$ to~$c$ becomes the rightmost geodesic~$\ph$ from~$h_0$ to~$c'$. The image of~$\ph'$ becomes~$\pr'$ in~$\tilde\m$. If $\ph\neq\rev(\pr')$, then it has to use one of the purple paths; the~$\le$ symbol (resp.~$<$ symbol) indicates that the length of the purple path is less than (resp.\ strictly less than) the length of the circumvented part. Tracking such purple paths back in~$\m$ shows that their existences contradict the definition of~$\ph'$.}
	\label{pftransfer}
\end{figure}

We prove that $\Phi_{\operatorname{right}}$ takes its values in~$\cM$ and that $\Phi_{\operatorname{left}}\circ\Phi_{\operatorname{right}}$ is the identity on~$\tilde\cM$ by the same argument.
\end{pre}
}

\subsection{Transferring from a face of degree one}

We now turn to\sv{ the setting of} Proposition~\ref{propp1m10}.\lv{ We let $\ba=(a_1,\ldots, a_r,1)$ be an $r+1$-tuple of positive integers with two odd coordinates and define $\tilde \ba=(\tilde a_1,\ldots,\tilde a_r)\de( a_1+1, a_2,\ldots, a_r)$.} We let~$\cM$ be the set of plane maps of type~$\ba$ carrying one distinguished corner in the first face and we let~$\tilde\cM$ be the set of plane maps of type~$\tilde\ba$ carrying one distinguished vertex and one distinguished half-edge incident to the first face and directed toward the distinguished vertex.\lv{ By Proposition~\ref{propbq}, the cardinalities of~$\cM$ and~$\tilde\cM$ are the sides of~\eqref{eqp1m10}.}

\sv{
\begin{center}
		\psfrag{f}[][][.8]{$f_1$}
		\psfrag{g}[][][.8]{$f_{r+1}$}
		\psfrag{1}[][][.7]{\textcolor{red}{$+1$}}
		\psfrag{v}[][][.8]{\textcolor{violet}{$v$}}
		\psfrag{M}[][][.8]{$\cM$}
		\psfrag{N}[][][.8]{$\tilde\cM$}
		\psfrag{c}[][][.8]{\textcolor{blue}{$c$}}
		\psfrag{h}[][][.8]{$\textcolor{blue}{h}\dt \textcolor{violet}{v}$}
	\centering\includegraphics[width=13cm]{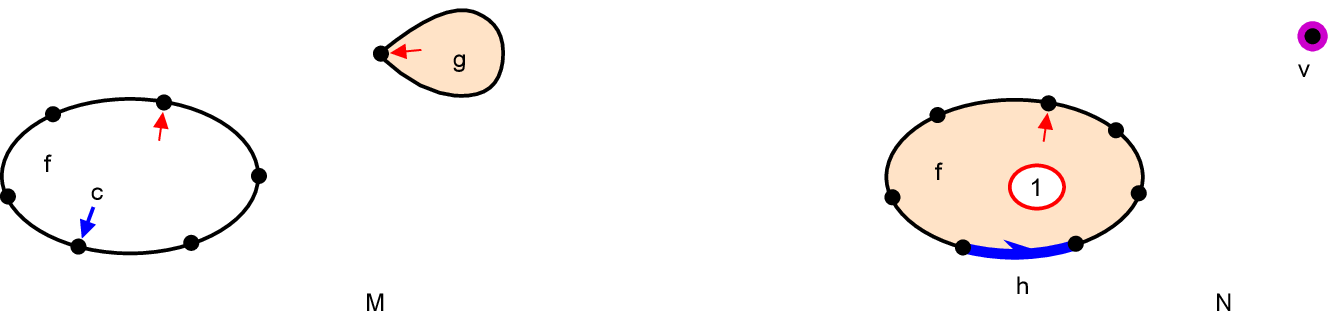}
\end{center}
}\lv{
\begin{center}
		\psfrag{f}[][]{$f_1$}
		\psfrag{g}[][]{$f_{r+1}$}
		\psfrag{1}[][][.9]{\textcolor{red}{$+1$}}
		\psfrag{v}[][][.8]{\textcolor{violet}{$v$}}
		\psfrag{M}[][]{$\cM$}
		\psfrag{N}[][]{$\tilde\cM$}
		\psfrag{c}[][]{\textcolor{blue}{$c$}}
		\psfrag{h}[][]{$\textcolor{blue}{h}\dt \textcolor{violet}{v}$}
	\centering\includegraphics[width=.95\linewidth]{pictoelem0}
\end{center}
}

The mappings are very similar as above; see\sv{ the right part of} Figure~\lv{\ref{transfer1}}\sv{\ref{transfer2}}. Let $(\m;c)\in \cM$. We slit~$\m$ along the rightmost geodesic~$\pp$ from the unique corner of~$f_{r+1}$ to~$c$. We denote by~$\pl$ and~$\pr$ the left and right copies of~$\pp$ in the resulting map and define~$\pr_0$ as the unique half-edge incident to~$f_{r+1}$. We then sew $\pl_{1\to [\pp]}$ onto $\pr_{0\to [\pp]-1}$, suppressing~$f_{r+1}$ in the process. In the resulting map, we denote by~$h$ the half-edge~$\rev(\pr)_{1}$ and denote by~$v$ the vertex~$\pl_{1}^-$. We then denote by~$\tilde\m$ the resulting map and let the outcome of the construction be $\Phi_{\operatorname{right}}^1(\m;c)\de(\tilde\m;v,h)$.

Conversely, starting from $(\tilde\m;v,h)\in \tilde\cM$, we consider the corner~$h_0$ delimited by~$h$ and its predecessor in the contour of~$f_1$, and denote by~$\pp'$ the leftmost geodesic from~$h_0$ to~$v$. We slit~$\tilde\m$ along~$\pp'$ starting from~$h_0$ and stopping at~$v$, \emph{without disconnecting the map} at~$v$. We denote by~$\pl'$ and~$\pr'$ the left and right copies of~$\pp'$ in the resulting map and sew $\pl'_{2\to [\pp']}$ onto $\pr'_{1\to [\pp']-1}$, thus creating a new degree~$1$-face, which we denote by~$f_{r+1}$. In the resulting map, we suppress the dangling edge~$\pl'_1$ and denote by~$c$ the corner it defines. We then denote by~$\m$ the resulting map and let the outcome of the construction be $\Phi_{\operatorname{left}}^1(\tilde\m;v,h)\de(\m;c)$.

\lv{
\begin{figure}[ht!]
		\psfrag{f}[][][.8]{$f_1$}
		\psfrag{g}[][][.8]{$f_{r+1}$}		
		\psfrag{p}[][][.8]{\textcolor{red}{$\pp'$}}
		\psfrag{l}[][][.8]{\textcolor{red}{$\pl'$}}
		\psfrag{r}[][][.8]{\textcolor{red}{$\pr'$}}
		
		\psfrag{m}[][][.8]{$\m$}
		\psfrag{n}[][][.8]{$\tilde\m$}
		
		\psfrag{c}[][][.8]{\textcolor{blue}{$c$}}
		\psfrag{h}[][][.8]{\textcolor{blue}{$h$}}
		\psfrag{e}[][][.8]{\textcolor{blue}{$h_0$}}
		\psfrag{v}[][][.8]{\textcolor{violet}{$v$}}
	\centering\includegraphics[width=8cm]{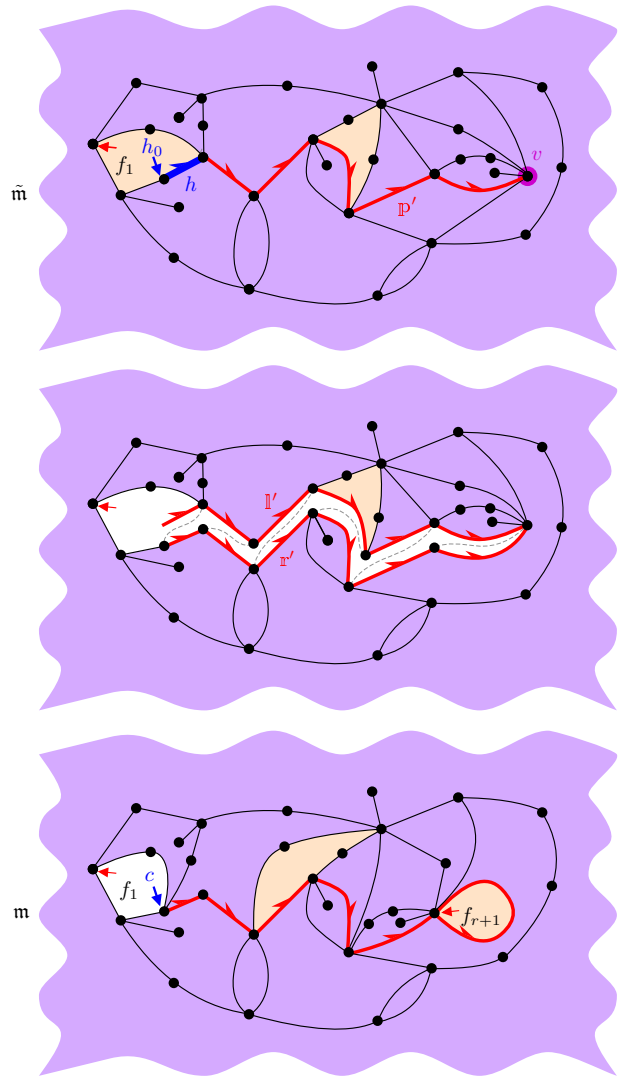}
	\caption{The transfer mappings in the case of a degree~$1$-face.}
	\label{transfer1}
\end{figure}
}

\begin{thm}
The mappings $\Phi_{\operatorname{right}}^1:\cM\to\tilde\cM$ and $\Phi_{\operatorname{left}}^1:\tilde\cM\to\cM$ are inverse bijections.
\end{thm}

\lv{
\begin{pre}
The proof is very similar to that of Theorem~\ref{thmtransfer2}; we leave it to the reader.
\end{pre}
}

\subsection{Decomposition of growing bijections into transfer bijections}\label{secdecomp}

\begin{figure}[ht]		
		\psfrag{m}[][][.8]{$\m$}
		\psfrag{n}[][][.8]{$\m'$}
		\psfrag{o}[][][.8]{$\m''$}
		\psfrag{p}[][][.8]{$\tilde\m$}
		\psfrag{f}[][][.5]{}
		
		\psfrag{c}[][][.8]{\textcolor{blue}{$c$}}
		\psfrag{d}[bl][bl][.8]{\textcolor{vert}{$c'$}}
		
		\psfrag{h}[][][.8]{\textcolor{blue}{$h$}}
		\psfrag{i}[][][.8]{\textcolor{vert}{$h'$}}
		\psfrag{e}[][][.8]{\textcolor{violet}{$e$}}
		\psfrag{v}[][][.8]{\textcolor{violet}{$v$}}
		\psfrag{j}[t][t][.8]{\textcolor{violet}{$h''$}}
	\centering\includegraphics[width=14cm]{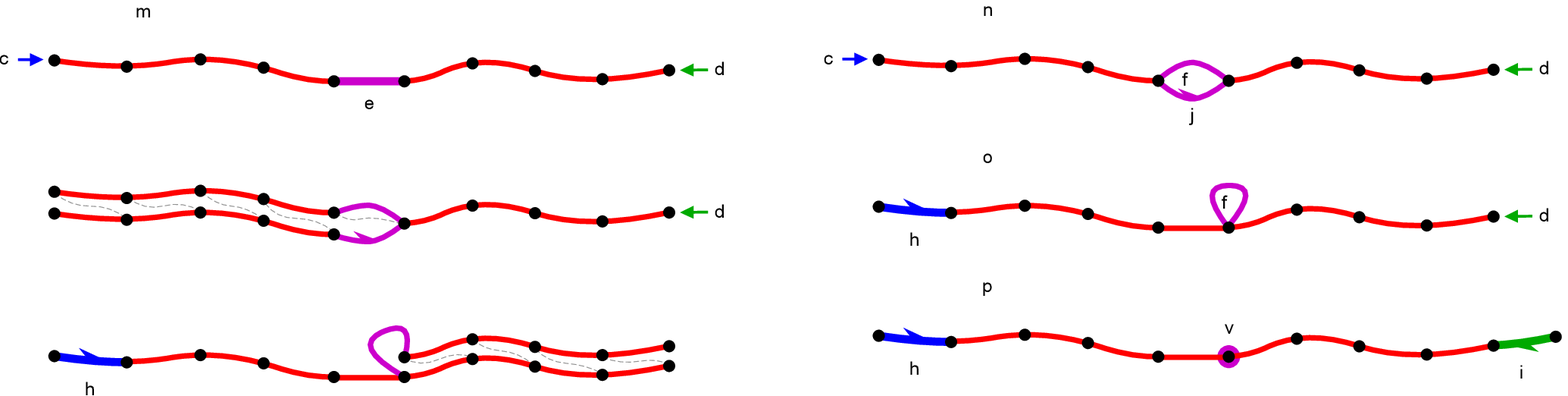}
	\caption{Two-step decomposition of a growing bijection into transfer bijections\lv{ in the simple case}.}
	\label{decomp}
\end{figure}

Let us explain our claim that growing bijections are compositions of two transfer bijections. We fix an $r$-tuple $\ba=(a_1,\ldots, a_r)$ of positive even integers and consider a map~$\m$ of type~$\ba$ with a distinguished edge~$e$ and two distinguished corners~$c$ and~$c'$ (either of the same face or of two different faces). We first define the map~$\m'$ of type $(a_1,\ldots, a_r,2)$ by replacing the distinguished edge~$e$ with an $r+1$-th face~$f_{r+1}$ of degree~$2$ by doubling the edge; the marked corner of this face is arbitrarily chosen. Next, we let~$h''$ be the unique half-edge incident to~$f_{r+1}$ that is directed away from~$c$. We set $(\m'';c'',h)\de\Phi_{\operatorname{right}}(\m';c,h'')$ and keep track of~$c'$ in the resulting map. The map~$\m''$ is of type $(a_1+1,a_2,\ldots, a_r,1)$ and we finally set $(\tilde\m;v,h')\de\Phi_{\operatorname{right}}^1(\m'';c')$, while keeping track of~$h$ in the resulting map. See Figure\lv{s}~\ref{decomp}\lv{ and~\ref{decomppinch}}.

\lv{
\begin{figure}[ht]		
		\psfrag{m}[][]{$\m$}
		\psfrag{n}[][]{$\m'$}
		\psfrag{o}[][]{$\m''$}
		\psfrag{p}[][]{$\tilde\m$}
		\psfrag{f}[][][.5]{}
		
		\psfrag{c}[][][.8]{\textcolor{blue}{$c$}}
		\psfrag{d}[bl][bl][.8]{\textcolor{vert}{$c'$}}
		
		\psfrag{h}[][][.8]{\textcolor{blue}{$h$}}
		\psfrag{i}[][][.8]{\textcolor{vert}{$h'$}}
		\psfrag{e}[][][.8]{\textcolor{violet}{$e$}}
		\psfrag{v}[][][.8]{\textcolor{violet}{$v$}}
		\psfrag{j}[l][l][.8]{\textcolor{violet}{$h''$}}
	\centering\includegraphics[width=15cm]{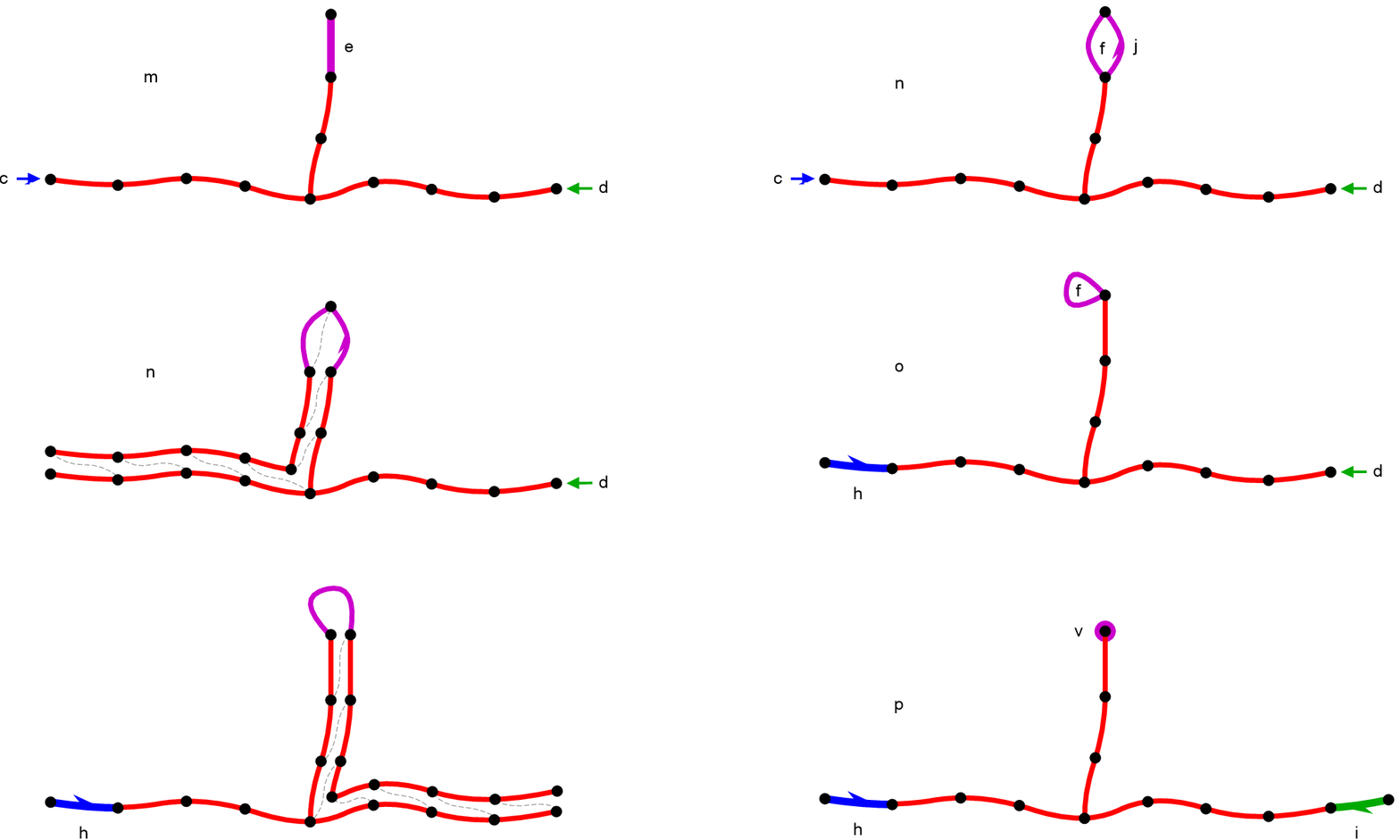}
	\caption{Two-step decomposition of a growing bijection into transfer bijections\lv{ in the pinched case}.}
	\label{decomppinch}
\end{figure}
}

We claim that $(\tilde\m;v,h,h')$ is exactly the output of the growing bijection of Section~\ref{secbb} or~\ref{secbq}. In~$\m$, the growing bijection uses two geodesics, one directed toward~$c$ and one directed toward~$c'$. Plainly, in the application of~$\Phi_{\operatorname{right}}$ to $(\m';c,h'')$, the sliding path in~$\m'$ corresponds to the geodesic directed toward~$c$. In order to show the claim, we only need to check that the image in~$\m''$ of the geodesic directed toward~$c'$ corresponds to the sliding path used by~$\Phi_{\operatorname{right}}^1$. This is because the mapping~$\Phi_{\operatorname{right}}$ only alters the map along the geodesic directed toward~$c$, which, by definition, cannot cross the geodesic directed toward~$c'$.

\lv{
\section{Uniform sampling}\label{secsample}

Our bijections can be used in order to sample a uniform bipartite or quasibipartite map of a given type~$\ba$. More precisely, let $\ba=(a_1,\ldots, a_r)$ be an $r$-tuple of positive even integers. Let~$\ba^1\de (2)$, $\ba^{n}\de \ba$ and~$\ba^2$, \ldots, $\ba^{n-1}$ be tuple of positive even integers such that, for all $2\le i\le n$, $\ba^i$ is obtained from~$\ba^{i-1}$
\begin{itemize}
	\item either by adding~$2$ to exactly one of its coordinates;
	\item or by concatenating it with~$(2)$.
\end{itemize}

In words, a map of type~$\ba^i$ differs from a map of type~$\ba^{i-1}$ by the fact that it has either one face having degree~$2$ more or one extra face of degree~$2$. For instance, one might choose the sequence
\begin{equation}\label{seqsample}
(2),\ (4),\ \ldots,\ (a_1),\ (a_1,2),\ (a_1,4),\ \ldots,\ (a_1,a_2),\ (a_1,a_2,2),\ \ldots,\ (a_1,\ldots, a_r)\,.
\end{equation}

We now sample a sequence of maps~$\m^1$, \ldots, $\m^n$ such that, for every $1\le i\le n$, the map~$\m^i$ is uniformly distributed among maps of type~$\ba^i$. Take for~$\m^1$ the only map of type~$\ba^1$. Then, sample~$\m^{i}$ from~$\m^{i-1}$ as follows.
\begin{itemize}
	\item If $a^{i}_j=a^{i-1}_j+2$, then choose uniformly at random in~$\m^{i-1}$ an edge and two corners in the $j$-th face and apply the mapping of Section~\ref{secp2+} (forgetting the distinguished elements in the resulting map).
	\item If $\ba^i$ is the concatenation of~$\ba^{i-1}$ with~$(2)$, then choose uniformly at random in~$\m^{i-1}$ an edge and transform it into a degree $2$-face, whose marked corner is uniformly chosen.
\end{itemize}

As the number of ways to choose the desired distinguished elements in a map only depends on the type of the map, distinguishing elements does not bias the uniform probability: if~$\m^{i-1}$ is uniformly distributed among the maps of type~$a^{i-1}$, then the map with uniformly chosen distinguished elements is uniformly distributed among maps of type~$a^{i-1}$ with the desired distinguished elements. As the mappings we use are bijections (either that of Section~\ref{secp2+} or the trivial one that changes a distinguished edge into an extra $2$-face), the resulting map with its distinguished elements is uniformly distributed among maps of type~$a^{i}$ with some distinguished elements. Finally, forgetting the distinguished elements yields a uniform map of type~$a^{i}$. By induction, for every $1\le i\le n$, the map~$\m^i$ is indeed uniformly distributed among maps of type~$\ba^i$, so that~$\m^n$ is uniformly distributed among maps of type~$\ba$, as desired.

Note that an advantage of choosing the sequence~\eqref{seqsample} is that we obtain a subsequence of ``growing'' uniform maps where the faces are added one by one. Namely, the map $\m^{a_1/2}$ is of type $(a_1)$, the map $\m^{a_1/2+a_2/2}$ is of type $(a_1,a_2)$, the map $\m^{a_1/2+a_2/2+a_3/2}$ is of type $(a_1,a_2,a_3)$, and so on. For instance, we may obtain in this way a sequence of uniform $2p$-angulations (maps of type $(2p,2p,\ldots,2p)$) of size $1$, $2$, $3$, \ldots, $n$ such that two subsequent maps do not differ too much.

Moreover, one can build on an already sampled uniform map in order to sample a larger one instead of starting from zero.

\bigskip

Now, in order to sample a uniform quasibipartite map, we proceed similarly, working with bipartite maps until the last step. Let~$\ba$ be a tuple of positive integers with two odd coordinates. Define~$\tilde\ba$ by adding one to an odd coordinate and subtracting one from the other odd coordinate (forget the null coordinate if there is one). Then sample from the algorithm above a uniform map of type~$\tilde\ba$ and use the mapping of Proposition~\ref{propp1m1} or~\ref{propp1m10} in order to obtain from it a uniform map of the desired type~$\ba$.

\section{Open questions and further discussion}\label{secopen}

In fact, the statements of Propositions~\ref{propp2}, \ref{propp1p1} and~\ref{propp1m1} are still valid as they are in the case of quasibipartite maps, that is, whenever~$\ba$ and~$\tilde\ba$ have at most two odd coordinates. We were not able to bijectively interpret this.

For Propositions~\ref{propp2} and~\ref{propp1p1} with a quasibipartite map on the left, that is, when~$\ba$ has two odd coordinates, the left-hand sides of~\eqref{eqp2} and~\eqref{eqp1p1} can still be interpreted as counting maps of type~$\ba$ carrying one distinguished edge and two distinguished corners. The distinguished edge may now be parallel to the distinguished corners and, in particular, it may very well be a loop. We do not see at the moment how to slit and slide when there is a loop on the sliding path.  

About Proposition~\ref{propp1m1}, when~$a_{r+1}$ is even and the map is quasibipartite, by Proposition~\ref{propbq}.\ref{pii}, two cases may happen. Either no half-edge incident to~$f_{r+1}$ is parallel to the distinguished corner, or exactly two half-edges incident to~$f_{r+1}$ are parallel to the distinguished corner. The term $\big\lfloor {a_{r+1}}/{2}\big\rfloor$ thus does not count half-edges incident to the last face
and directed toward a distinguished corner~$c$. One might need to add one of the two parallel half-edges in this case, for instance, the one parallel half-edge~$h$ such that~$c$ lies to the right of the loop made up by the two rightmost geodesics from~$h$ and from~$\rev(h)$ to~$c$, oriented by~$h$. 

Another difficulty is foreseeable in the setting of Proposition~\ref{propp2} when~$a_1$ is even and two coordinates of~$\ba$ are odd. If we hope to find a slit-slide-sew bijection that can be decomposed as two transfer bijections, one will need to exit the realm of bipartite or quasibipartite maps when transferring a corner from the extra degree-2 face to the first face\ldots
}

\bibliographystyle{alpha}
\bibliography{main}

\end{document}